\documentclass{amsart}

\usepackage{amsmath}
\usepackage{amsfonts}
\usepackage{amssymb,enumerate}
\usepackage{amsthm}
\usepackage[all]{xy}

\newtheorem{lem}{Lemma}[section]
\newtheorem{cor}[lem]{Corollary}
\newtheorem{prop}[lem]{Proposition}
\newtheorem{thm}[lem]{Theorem}

\newtheorem{Defn}[lem]{Definition}
\newtheorem{Ex}[lem]{Example}
\newtheorem{Question}[lem]{Question}
\newtheorem{Property}[lem]{Property}
\newtheorem{Properties}[lem]{Properties}
\newtheorem{Discussion}[lem]{Remark}
\newtheorem{Construction}[lem]{Construction}
\newtheorem{Subprops}{}[lem]
\newtheorem{Para}[lem]{}

\newenvironment{question}{\begin{Question}\rm}{\end{Question}}

\newenvironment{para}{\begin{Para}\rm}{\end{Para}}

\newcommand{\y}{\mathbf{y}}
\newcommand{\cbc}[2]{#1(#2)}
\newcommand{\wt}{\widetilde}

\newcommand{\comp}[1]{\widehat{#1}}
\newcommand{\ideal}[1]{\mathfrak{#1}}
\newcommand{\m}{\ideal{m}}
\newcommand{\n}{\ideal{n}}
\newcommand{\p}{\ideal{p}}
\newcommand{\q}{\ideal{q}}
\newcommand{\fa}{\ideal{a}}
\newcommand{\fb}{\ideal{b}}
\newcommand{\fN}{\ideal{N}}

\newcommand{\fr}{\ideal{r}}
	
\newcommand{\pd}{\mathrm{pd}}	
	
\newcommand{\gkdim}[1]{\mathrm{G}_{#1}\text{-}\dim}	
\newcommand{\ext}{\mathrm{Ext}}	
\newcommand{\Ht}{\mathrm{ht}}	
	
\newcommand{\depth}{\mathrm{depth}}	
\newcommand{\rank}{\mathrm{rank}}	
\newcommand{\zz}{\mathbb{Z}}

\newcommand{\id}{\mathrm{id}}	
\newcommand{\rhom}{\mathbf{R}\mathrm{Hom}}	
\newcommand{\lotimes}{\otimes^{\mathbf{L}}}
\newcommand{\amp}{\mathrm{amp}}
\newcommand{\HH}{\mathrm{H}}
\newcommand{\Hom}{\mathrm{Hom}}

\newcommand{\fd}{\mathrm{fd}}

\newcommand{\spec}{\mathrm{Spec}}
\newcommand{\mspec}{\mathrm{m\text{-}Spec}}

\newcommand{\s}{\mathfrak{S}}
\newcommand{\x}{\mathbf{x}}
\newcommand{\tor}{\mathrm{Tor}}
\newcommand{\vf}{\varphi}

\newcommand{\D}{\mathsf{D}}

\newcommand{\len}{\operatorname{length}}

\newcommand{\xra}{\xrightarrow}
\newcommand{\shift}{\mathsf{\Sigma}}

\newcommand{\Cl}{\operatorname{Cl}}

\newcommand{\ass}{\operatorname{Ass}}
\newcommand{\Pic}{\operatorname{Pic}}

\newcommand{\card}{\operatorname{card}}

\newcommand{\tri}{\trianglelefteq}
\newcommand{\Mod}{\operatorname{Mod}}
\newcommand{\Ker}{\operatorname{Ker}}
\newcommand{\minh}{\operatorname{Minh}}
\newcommand{\Min}{\operatorname{Min}}
\newcommand{\Div}{\operatorname{Div}}
\newcommand{\divmap}{\operatorname{div}}
\newcommand{\Prin}{\operatorname{Prin}}

\begin{document}

\bibliographystyle{amsplain}

\author{Sean Sather-Wagstaff}

\address{Sean Sather-Wagstaff, Department of Mathematics, California
  State University, Dominguez Hills, 1000 E.~Victoria St., Carson, CA
  90747 U.S.A.}

\curraddr{Department of Mathematical Sciences, Kent State University,
  Mathematics and Computer Science Building, Summit Street, Kent OH
  44242}

\email{sather@math.kent.edu}

\urladdr{http://www.math.kent.edu/~sather}

\thanks{This research 
was conducted in part while the author was an NSF Mathematical Sciences 
Postdoctoral Research Fellow and a visitor at the University of 
Nebraska-Lincoln.}

\title{Semidualizing modules and the divisor class group}

\dedicatory{Dedicated to Phillip Griffith on the occasion of his retirement}

\keywords{semidualizing complexes, semidualizing modules, G-dimensions,
Gorenstein dimensions, determinantal rings, divisor class groups}
\subjclass[2000]{Primary 13C05, 13C13, 13C20; Secondary 13C40, 13D05, 
13D25}

\begin{abstract}
Among the finitely generated modules over a Noetherian ring $R$, the
semidualizing modules have been singled out due to their
particularly nice duality properties.  
When $R$ is a normal domain, we exhibit a natural inclusion of the set of 
isomorphism classes of semidualizing $R$-modules into the divisor class 
group of $R$.  After a description of the basic properties of this 
inclusion, it is employed to
investigate the 
structure of 
the set of 
isomorphism classes of semidualizing $R$-modules.
In particular, 
this set is described completely
for determinantal
rings over normal domains.
\end{abstract}

\maketitle

\section{Introduction} 

Semidualizing modules arise naturally in the investigations
of various 
duality theories in commutative algebra.  One instance of this 
is Grothendieck and Hartshorne's local duality wherein a dualizing module,
or more generally a dualizing complex, is employed to study local 
cohomology~\cite{hartshorne:rad,hartshorne:lc}.  Another instance is 
Auslander and Bridger's methodical study of duality properties with 
respect to a rank 1 free module that gives rise to the Gorenstein 
dimension~\cite{auslander:adgeteac,auslander:smt}.  
A free module of rank 1 and 
a dualizing module are both examples of semidualizing modules.

Let $R$ be a Noetherian ring.  A finitely generated $R$-module $C$ is 
\emph{semidualizing} if the natural homothety map 
$R\to\Hom_R(C,C)$ is an isomorphism and $\ext^i_R(C,C)=0$ for 
each integer $i>0$.  The study of such modules in the abstract was 
initiated by Foxby~\cite{foxby:gmarm} and Golod~\cite{golod:gdagpi} 
where they were called ``suitable'' modules, and has been continued 
recently by others; see for 
example~\cite{christensen:scatac, 
frankild:sdcms,gerko:sdc,gerko:osmagi}.  

The semidualizing modules 
and more generally the semidualizing \emph{complexes} 
are useful, for example, in identifying local homomorphisms of 
finite Gorenstein dimension with particularly nice 
properties
as in~\cite{avramov:rhafgd,frankild:qcmpolh}.
This utility, along with our desire to expand upon it,
motivates our investigation of the basic properties of 
such modules and  the structure of the 
entire set of 
isomorphism classes of semidualizing $R$-modules, which we denote
$\s_0(R)$.  Surprisingly little is known about his set.
For instance, 
researchers in this subject have been grappling with the following open question 
for several years now;  see~\cite[(1.2)]{christensen:dscrap} for recent progress.

\begin{question} \label{q01}
When $R$ is local, must $\s_0(R)$ be finite?
\end{question}

This work 
is part of a research effort focused on determining the overall structure of $\s_0(R)$.
Much of the ground work for this effort, motivated by~\cite{christensen:scatac, gerko:sdc},
is found in~\cite{frankild:appx}.  
Initial evidence of the richness of the structure of this set is found in 
the fact that it admits an ordering described in terms of a reflexivity relation;  
see~\ref{paranew4}.  
Further structure is uncovered in~\cite{frankild:sdcms} where
numerical data from this ordering 
is used to build a nontrivial metric on $\s_0(R)$.
While the existence of a metric does not itself provide answers to any
of the open questions about the structure of $\s_0(R)$, it represents a new perspective from which 
to view this set.  This perspective has proved particularly useful 
for identifying
questions that we would not have otherwise thought to ask.
For instance, our investigation into the nontriviality of the metric
led us to the fact that, when $R$ is Cohen-Macaulay and $\s_0(R)$ is nontrivial,
there exist elements of $\s_0(R)$ that are incomparable under the reflexivity ordering;
see~\cite[(3.5)]{frankild:sdcms}.

In the current paper, we forward another new perspective from which to 
investigate the set $\s_0(R)$.
It is motivated by Bruns' work~\cite{bruns:1} wherein the divisor class group is used
to describe the dualizing module of certain Cohen-Macaulay normal domains.
Accordingly, when $R$ is a normal domain,
we exhibit a natural inclusion
$\s_0(R)\hookrightarrow\Cl(R)$ that behaves well with respect to 
standard operations.  
This inclusion allows us to exploit the 
known behavior of $\Cl(R)$ to gain insight into the 
structure of $\s_0(R)$.  For instance, 
if $\Cl(R)$ is finite, then $\s_0(R)$ is also finite.  
The basic results from this analysis are presented in Section~\ref{sec2}.  

Section~\ref{sec3} contains the meat of this investigation and demonstrates
the power of this new perspective. 
It consists of
analyses showing how the divisor class group can be used to 
give a complete description of the set of semidualizing modules for certain classes of rings.
We recount here three such situations, focusing our attention on the 
finiteness and size of $\s_0(R)$.
First, taking our cues 
from~\cite{bruns:1}, we describe the semidualizing modules over a 
determinantal ring in Theorem~\ref{thm03}.  

\begin{thm} \label{thmA}
Let $A$ be a 
normal domain and 
$m,n,r$ nonnegative integers such that
$r<\min\{m,n\}$.  With
$\mathbf{X}=\{X_{ij}\}$ an $m\times n$ matrix of
variables, set 
$R=A[\mathbf{X}]/I_{r+1}(\mathbf{X})$
where $I_{r+1}(\mathbf{X})$
is the ideal generated by the minors of $\mathbf{X}$ of size $r+1$.  
The set $\s_0(R)$ is finite if and 
only if $\s_0(A)$ is so.  More specifically, one has the following cases.
\begin{enumerate}[\quad\rm(a)]
\item 
If $r=0$ or $m=n$, then there is a bijection $\s_0(R)\approx\s_0(A)$.
\item 
If $r>0$ and $m\neq n$, then there is a bijection
$\s_0(R) \approx\s_0(A)\times\{0,1\}$.
\end{enumerate}
\end{thm}

When $A$ is 
a graded Cohen-Macaulay (super-)normal domain with $A_0$ local 
(and complete), 
this result extends to the localization (and to the completion) of
$R$ at its graded maximal ideal;  see Corollaries~\ref{cor05} 
and~\ref{cor06}.  

The next two results demonstrate how 
this technique yields information about rings that are themselves 
not normal domains;
they are contained in
Corollary~\ref{cor09}.  
Theorem~\ref{thmC} is new even when $B$ is a field.

\begin{thm} \label{thmB}
Let $A=\coprod_{i\geq 0}A_i$ be a graded super-normal 
domain 
with $A_0$ local and complete.  Let $\n$ be the graded maximal ideal 
of $A$ and $\comp{A}$ the $\n$-adic 
completion
of $A$.  Let $\mathbf{y}=y_1,\ldots,y_q\in\n A_{\n}$ be an 
$A_{\n}$-sequence and fix an integer $m\geq 1$.  
There are bijections
$$\s_0(\comp{A}/(\mathbf{y})^m)\approx\s_0(A_{\n}/(\mathbf{y})^m)\approx
\begin{cases}
\s_0(A_{\n}) & \text{if $m=1$ or $q=1$} \\
\s_0(A_{\n})\times\{0,1\} & \text{if $m,q>1$.}
\end{cases}$$
\end{thm}

\begin{thm}\label{thmC}
With $A$ as in Theorem~\ref{thmB}, let $B$ denote either 
$A_{\n}$ or $\comp{A}$,
and let $t$ be a positive integer.
For $l=1,\ldots,t$ fix a positive integer $q_l$ and
set
\[ 
S=(B\ltimes B^{q_1})\otimes_{B}\cdots\otimes_{B}(B\ltimes 
B^{q_t}).
\]
If $s$ is the number of indices $l$ with $q_l> 1$,
then and there is a bijection
\[ \s_0(B)\times\{0,1\}^s 
\approx
\s_0(S). \]
\end{thm}

The statements of our main results are module-theoretic in nature.  
However, we often employ tools from the derived category.  
We include a summary of the relevant notions in
Section~\ref{sec1} along with basic facts about semidualizing modules 
and the divisor class group. 

\section{Background} \label{sec1}

In this paper, the term ``ring'' is used for a commutative Noetherian ring with
identity, and ``module'' is used for a unital module.  Let $R$ be a ring.

\begin{para} \label{paranew1}
An \emph{$R$-complex} is a sequence of $R$-module homomorphisms
\[ X=\cdots\xra{\partial^X_{i+1}}X_i\xra{\partial^X_{i}}X_{i-1} 
\xra{\partial^X_{i-1}}\cdots \]
with $\partial_i^X\partial_{i+1}^X=0$ for each $i$.  We work occasionally
in the derived category $\D(R)$ whose objects are the $R$-complexes;   
references on the subject
include~\cite{gelfand:moha,hartshorne:rad,verdier:cd,verdier:1}.
The category of $R$-modules $\Mod(R)$ 
is naturally identified with the full 
subcategory of $\D(R)$ whose objects are the complexes 
homologically concentrated in 
degree 0.
For $R$-complexes $X$ and $Y$ 
the left derived tensor product complex 
is denoted $X\lotimes_R Y$ and 
the right derived homomorphism complex is 
$\rhom_R(X,Y)$.  For an integer $n$, 
the $n$th \emph{shift} or \emph{suspension} of $X$ is denoted
$\shift^n X$ where $(\shift^n 
X)_i=X_{i-n}$ and $\partial_i^{\shift^n X}=(-1)^n\partial_{i-n}^X$.
The symbol ``$\simeq$'' indicates an
isomorphism in $\D(R)$, and ``$\sim$'' is isomorphism up to shift.

A complex $X$ is \emph{homologically finite}, respectively 
\emph{homologically degreewise finite}, if its total 
homology module $\HH(X)$, respectively each individual homology 
module $\HH_i(X)$, is a finite $R$-module.
The \emph{infimum}, \emph{supremum}, and \emph{amplitude} of 
$X$ are 
\begin{align*}
\inf(X)&=\inf\{i\in\zz\mid\HH_i(X)\neq 0\} \\
\sup(X)&=\sup\{i\in\zz\mid\HH_i(X)\neq 0\} \\
\amp(X)&=\sup(X)-\inf(X) 
\end{align*}
respectively, with the conventions $\inf\emptyset = \infty$ and
$\sup\emptyset = -\infty$.
When $R$ is local with residue field $k$,
the \emph{depth} of a homologically finite complex $X$ is
$$\depth_RX=-\sup(\rhom_R(k,X)).$$
The \emph{Bass series} of $X$ is the formal 
Laurent series
$I_R^X(t)=\sum_i\mu^i_R(X)t^i$ where
\[
\mu_R^i(X)=\rank_k\HH_{-i}(\rhom_R(k,X)) \]
for each integer $i$.   From Foxby~\cite[(13.11)]{foxby:hacr} the quantity $\id_R 
X$ is finite if and only if $I_R^X(t)$ is a Laurent polynomial.  
\end{para}

\begin{para} \label{paranew2}
Associated to a complex $K$ is a
natural homothety morphism
\[ \chi^R_K  \colon R\to\rhom_R(K,K). \]
When $K$ is homologically finite, it is \emph{semidualizing} if 
$\chi^R_K$ is an isomorphism.
A complex $D$ is \emph{dualizing} if it is semidualizing and has 
finite injective dimension.
The set of shift-isomorphism classes of semidualizing $R$-complexes 
is denoted $\s(R)$, and the class of a semidualizing 
complex $K$ in $\s(R)$
is denoted $[K]_R$ or simply $[K]$ when there is no 
danger of confusion.  The ring $R$ is \emph{$\s$-finite} 
if $\s(R)$ is a finite set.

A finitely generated $R$-module $C$ is 
\emph{semidualizing} if the natural homothety map 
$R\to\Hom_R(C,C)$ is an isomorphism and $\ext^{\geqslant 1}_R(C,C)=0$.
The module $R$ is 
semidualizing.  When $R$ is Cohen-Macaulay, a \emph{dualizing module} 
(or \emph{canonical module})
is a semidualizing module of finite injective dimension.
The set of isomorphism classes of semidualizing $R$-modules 
is denoted $\s_0(R)$.  The identification of $\Mod(R)$
with a 
subcategory of $\D(R)$ provides a natural inclusion 
$\s_0(R)\hookrightarrow\s(R)$, and we  identify $\s_0(R)$ 
with its image in $\s(R)$.  In particular,
the class of a semidualizing 
module $C$ in $\s_0(R)$ is denoted $[C]_R$ or $[C]$.  
The ring $R$ is \emph{$\s_0$-finite} 
if $\s_0(R)$ is a finite set.

Some of our favorite ring theoretic properties have characterizations 
in terms of semidualizing objects.
If $R$ is Cohen-Macaulay local, then $\s(R)=\s_0(R)$.  
If $R$ is Gorenstein local, then $\s(R)=\{[R]\}$.  The converses 
hold when $R$ admits a dualizing complex;  
see Christensen~\cite[(3.7),(8.6)]{christensen:scatac}.
\end{para}

\begin{para} \label{paranew3}
Let $K$ be a semidualizing complex.  A homologically finite complex $X$ 
is \emph{$K$-reflexive} 
if
$\rhom_R(X,K)$ is homologically bounded and the natural
biduality morphism
\[ \delta^K_X  \colon X\to\rhom_R(\rhom_R(X,K),K) \]
is an isomorphism.  For instance, the complexes $R$ 
and $K$ are both $K$-reflexive.  When $\dim(R)$ is finite,
the complex $K$ is dualizing if and only if 
every homologically finite $R$-complex is $K$-reflexive.
The \emph{$\mathrm{G}_K$-dimension} of $X$ is 
\[ \gkdim{K}_R X = 
\begin{cases} 
\inf K - \inf\rhom_R(X,K) & \text{when $X$ is $K$-reflexive} \\
\infty & \text{otherwise.} \end{cases} \]
If $R$ is local and $X$ is $K$-reflexive,
then the AB 
formula~\cite[(3.14)]{christensen:scatac} reads
\[ \gkdim{K}_R X =\depth R -  \depth_R X.\]
When $C$ is a semidualizing module, the 
$\mathrm{G}_C$-dimension of a finitely generated $R$-module $M$
can be described in terms of resolutions.  
We first describe the modules used in the resolutions.
A finitely generated $R$-module $G$ is \emph{totally $C$-reflexive} if the 
natural biduality map $G\to\Hom_R(\Hom_R(G,C),C)$ is bijective, and 
$\ext_R^{\geqslant 1}(G,C)=0=\ext^{\geqslant 1}_R(\Hom_R(G,C),C)$.  A finitely generated 
$R$-module $M$ then has finite $\mathrm{G}_C$-dimension if and only 
if it admits a resolution 
\[ 0\to G_g\to\cdots\to G_0\to M\to 0 \]
with each $G_i$ totally $C$-reflexive;  the $\mathrm{G}_C$-dimension 
of $M$ is then the minimum integer $g$ for which $M$ admits such a resolution.
\end{para}

\begin{para} \label{paranew4}
The above notion of reflexivity gives rise to
orderings on the sets $\s_0(R)$ and $\s(R)$:  write
$[C]\tri[C']$ whenever $C'$ is $C$-reflexive.  
This ordering is trivially reflexive:  $[C]\tri[C]$.
Also, when $R$ is local,
the ordering
is antisymmetric:  if $[C]\tri[C']$ and $[C']\tri[C]$, then
$[C]=[C']$;  see~\cite[(5.3)]{takahashi:hiatsb}.  The question of the transitivity of this ordering
has been of interest in this area for some time:
\end{para}

\begin{question} \label{q05}
If $[C]\tri[C']$ and $[C']\tri[C'']$, then must one have
$[C]\tri[C'']$?
\end{question}

\begin{para} \label{paranew5}
For $1,\ldots,n$ 
let $U_i$ be a set with a relation
$\tri_i$.  The Cartesian product $U_1\times\cdots\times U_n$ is 
endowed with the
\emph{product relation}:  $(u_1,\ldots,u_n)\tri(u_1',\ldots,u_n')$  when
$u_i\tri_iu_i'$ for each $i=1,\ldots,n$.  
It follows immediately from the definition that
$\tri$ is transitive if and only if each $\tri_i$ is transitive. 
A map 
$\alpha\colon U_1\to U_2$ is \emph{order-respecting} when
$u\tri_1 u'$ implies $\alpha(u)\tri_2\alpha(u')$, and it is 
\emph{perfectly order-respecting} when the converse also holds.
Observe that, when $\alpha$ is a perfectly order-respecting
bijection, the relations $\tri_1,\tri_2$ are
simultaneously transitive.
The symbol $\approx$ indicates a perfectly order-respecting bijection.
\end{para}

We pose one final question in this section, motivated by the following well-known
equality of Hilbert-Samuel multiplicities: 
If $R$ is a local Cohen-Macaulay ring with dualizing module
$\omega$, then $e(\omega)=e(R)$.

\begin{question} \label{q06}
Let $(R,\m)$ be a local ring and $I$ an $\m$-primary ideal.  
If $C$ is a semidualizing $R$-module, must there be an equality
$e(I,C)=e(I,R)$?
\end{question}

Consult Matsumura~\cite[\S 14]{matsumura:crt} 
for the basics of the Hilbert-Samuel
multiplicity.  We answer Question~\ref{q06} in the affirmative
for several classes of rings in Corollaries~\ref{cor08}
and~\ref{cor09}.
To do so, we need the following lemma, which addresses this question
when $R$ is generically Gorenstein.

\begin{lem}
\label{lem:mult}
Let $(R,\m)$ be a local ring and $C$ a semidualizing $R$-module.
\begin{enumerate}[\quad\rm(a)]
\item 
\label{item:mult:a}
If $R_{\p}$ is Gorenstein for each $\p\in\spec(R)$ with
$\dim(R/\p)=\dim(R)$, then $e(J,C)=e(J,R)$ for each $\m$-primary ideal $J$.
\item
\label{item:mult:b}
Assume that $R$ is equidimensional and, for each $\p\in\Min(R)$, the rings
$R_{\p}$ and $R_{\p}/\p R_{\p}\otimes_R \comp{R}$ are Gorenstein.
If $\y\in\m$ is an $R$-sequence and $R'=R/\y$,
then $e(J,C\otimes R')=e(J,R')$ for each $\m R'$-primary ideal $J$.
\end{enumerate}
\end{lem}

\begin{proof}
\eqref{item:mult:a}
Set $\minh(R)=\{\p\in\spec(R)\mid\dim(R/\p)=\dim(R)$.  
For each $\p\in\minh(R)$, there is an isomorphism
$C_{\p}\cong R_{\p}$ by~\cite[(8.6)]{christensen:scatac} since $R_{\p}$ is Gorenstein.
This provides the second equality in the following sequence
$$e(J,C)=\hspace{-3mm}\sum_{\p\in\minh(R)}\hspace{-3mm}\len(C_{\p})e(J,R/\p)
=\hspace{-3mm}\sum_{\p\in\minh(R)}\hspace{-3mm}\len(R_{\p})e(J,R/\p)=e(J,R) $$
while the others are the additivity formulas for multiplicies~\cite[(4.7.t)]{bruns:cmr}.

\eqref{item:mult:b}
Recall the following fact:
If $(A,\m_A)\to (B,\m_B)$ is a flat local homomorphism
such that $\m_B=\m_AB$, and if $M$ is a finite $A$ module,
then for each $\m_A$-primary ideal $I$ there is an equality of Hilbert functions
$$\len_A(M/I^nM)=\len_B((M\otimes_A B)/(IB)^n(M\otimes_A B))$$
which yields an equality of multiplicities; see~\cite[(2.3)]{herzog:mlr}
\begin{equation} \tag{\ref{lem:mult}.1}
e(I,M)=e(IB,M\otimes_A B). \label{eq:mult} 
\end{equation}
Next, from~\cite[(0.10.3.1)]{grothendieck:ega3-1} one has a flat
local homomorphism $\vf\colon (R,\m)\to (S,\n)$ such that $S$ 
has infinite residue field
and $\n=\m S$.
Since $S$ is flat over $R$, the sequence $\y$ is $S$-regular.
Set $S'=S/(\y)S$ and let $\tau\colon S\to S'$ denote the natural surjection.
Note that 
the induced map $R'\to S'$ is flat and local and that
the maximal ideal of $R'$ extends to that of $S'$.
It follows that the ideal $JS'$ is $\m S'$-primary.

For every $\q\in\minh(S)$, 
the local ring
$S_{\q}$ is Gorenstein.  
To see this, fix  $\q\in\minh(S)$ and set $\p=\q\cap R$. 
Since $\vf$ is faithfully flat, the going-down theorem  implies
$\Ht(\p)\leq\Ht(\q)=0$, and so  
$\p\in\Min(R)$.  By assumption, the rings  $R_{\p}/\p R_{\p}\otimes_R \comp{R}$
and $S/\m S$
are Gorenstein, and it follows from~\cite[Main Thm.]{avramov:glp} 
that the fibre $R_{\p}/\p R_{\p}\otimes_R S$ is 
Gorenstein.  
Thus, the induced map $\vf_{\q}\colon R_{\p}\to S_{\q}$ is flat and local with Gorenstein
source and Gorenstein closed fibre.   This implies that $S_{\q}$ is Gorenstein.

The ring $S'$ has a  parameter ideal
$\mathbf{x}=(x_1,\ldots,x_r)S'\subseteq JS'$ such that 
\begin{equation} 
\label{eq:mult2} \tag{\ref{lem:mult}.2}
e(JS',C\otimes_{R}S')=e(\x,C\otimes_{R}S') \qquad\text{and}\qquad 
e(JS',S')=e(\x,S')
\end{equation}
see~\cite[(4.6.10)]{bruns:cmr}.
Fix $\wt{x}_1,\ldots,\wt{x}_r\in S$ such that
$\tau(\wt{x}_i)=x_i$ for each $i=1,\ldots,r$, and set $\wt{J}=(\wt{\x},\y)S$.
The third equality in the next sequence is from part~\eqref{item:mult:a}
\begin{equation*} 
e(\x,C\otimes_{R}S')
=e(\x,(C\otimes_{R}S)\otimes_SS')=e(\wt{J},C\otimes_{R}S)
=e(\wt{J},S)=e(\x,S').
\end{equation*}
The first one comes from the isomorphism 
$(C\otimes_{R}S)\otimes_SS'\cong C\otimes_{R}S'$,
and the
second and fourth hold by~\cite[(14.11)]{matsumura:crt}.
This yields the third equality below 
\begin{align*}
e(J,C\otimes_R R')
&=e(JS',C\otimes_{R}S')
=e(\x,C\otimes_{R}S')\\
&=e(\x,S')
=e(JS',S')
=e(J,R')
\end{align*}
while the remaining ones are from equations~\eqref{eq:mult} and~\eqref{eq:mult2}.
\end{proof}

\begin{para} \label{paranew6}
Let $\vf\colon R\to S$ be a ring homomorphism.
The \emph{flat dimension} of $\vf$ is $\fd(\vf)=\fd_R(S)$.
Assume that $\vf$ is surjective and
$\fd(\vf)<\infty$.  
The map $\vf$ is 
\emph{Cohen-Macaulay} of grade $d$ if
$S$ is a perfect $R$-module of grade $d$.
The map is 
\emph{Gorenstein} of grade $d$ 
if it is
Cohen-Macaulay of grade $d$ and, 
for each prime ideal $\q\subset S$, the $S_{\q}$-module
$\ext^{d}_R(S,R)_{\q}$ is cyclic.
\end{para}

\begin{para} \label{paranew7}
This section concludes with the definition of the divisor class group 
of a
normal domain $R$ with field of fractions $Q$.  Let $(-)^\star$ denote the functor 
$\Hom_R(-,R)$.
An $R$-module $M$ is \emph{reflexive}\footnote{Not to be confused with
``$K$-reflexive'' or ``totally $C$-reflexive''.}
if it is finitely generated and the natural biduality map 
$b_M^S\colon M\to M^{\star\star}$ is bijective.
The \emph{divisor class group} of $R$,
denoted $\Cl(R)$, is the 
set of isomorphism classes of reflexive $R$-modules of rank 1.
The 
isomorphism class of a 
reflexive module $M$ is denoted $[M]_R$ or $[M]$ when there is no 
risk of confusion.
The set $\Cl(R)$ admits an Abelian group structure:
when $M,N$ are rank 1 reflexive modules
\[ [M]+[N]=[(M\otimes_R N)^{\star\star}]
\qquad
[M]-[N]=[\Hom_R(N,M)]. \]
If $\fa,\fb$ are ideals with $\fa\cong M$ and $\fb\cong N$, then
$[M]+[N]=[\fa]+[\fb]=[(\fa\fb)^{\star\star}]$.  

The fact that the operations described above make $\Cl(R)$ into an Abelian group
seems to be part of the folklore of this subject; see~\cite[Sec.\ 0]{lipman}.
We sketch a proof of this fact below which also indicates why
this definition is equivalent to other formulations that 
may be more familiar to some readers.

A \emph{fractionary} ideal of $R$ is a nonzero finitely generated $R$-submodule of $Q$.
From the proof of~\cite[(2.2.iv)]{fossum:dcgkd}, one has $\Hom_R(\fa,\fb)\cong \fa:_Q\fb$
for every pair of fractionary ideals $\fa,\fb$.
Let $D(R)$ denote the set of reflexive fractionary ideals of $R$, and
let $P(R)$ denote the set of principal fractionary ideals of $R$.  
As $R$ is a normal domain, 
we learn from~\cite[(3.4)]{fossum:dcgkd} that
$D(R)$ is an
Abelian group via the operations
\[ \fa+\fb=R:_Q(R:_Q(\fa\fb))
\qquad
\fa-\fb=\fa:_Q \fb \]
with identity $R$.
One checks
that $P(R)$ is a subgroup of $D(R)$ and that $\overline{\fa}=\overline{\fb}$
in $D(R)/P(R)$ if and only if $\fa=a\fb$ for some $a\in Q^{\times}$ if and only if
$\fa\cong \fb$.

Let $Z(R)$ denote the set of all height 1 prime ideals $\p\subset R$.
For each $\p\in Z(R)$, the localization $R_{\p}$ is a discrete valuation ring
because $R$ is a normal domain, and we let $v_{\p}\colon Q^{\times}\to\mathbb{Z}$
denote the associated valuation.
For each $\fa\in D(R)$ and $\p\in Z(R)$, set
$v_{\p}(\fa)=\inf\{v_{\p}(a)\mid a\in\fa\}$.
Set $\Div(R)=\oplus_{\p\in Z(R)}\mathbb{Z}\cdot[R/\p]$
and consider the function  $\divmap\colon D(R)\to \Div(R)$
given by $\divmap(\fa)=(v_{\p}(\fa)[R/\p])_{\p\in Z(R)}$.
This is an Abelian group isomorphism
by~\cite[(5.9)]{fossum:dcgkd}
and we set  $\Prin(R)=\divmap(P(R))\subseteq \Div(R)$.
It is routine to verify that $\divmap$ induces a group isomorphism
$D(R)/P(R)\cong \Div(R)/\Prin(R)$.  

Many readers will undoubtedly recognize $\Div(R)/\Prin(R)$ as the 
definition of the divisor class group from~\cite{bourbaki:acc7d}.
To see that this is equivalent to the definition formulated above
(and that our formulation yields an Abelian group) it suffices
to construct a bijection $f\colon D(R)/P(R)\to \Cl(R)$ such that
$$f(\overline{\fa}+\overline{\fb})=f(\overline{\fa})+f(\overline{\fb})
\qquad f(\overline{\fa}-\overline{\fb})=f(\overline{\fa})-f(\overline{\fb})
\qquad f(\overline{R})=[R]$$
for each $\fa,\fb\in D(R)$.
Since $\overline{\fa}=\overline{\fb}$
if and only if $\fa\cong \fb$,
one sees that the assignment $\overline{\fa}\mapsto[\fb]$
describes a well-defined injection $f\colon D(R)/P(R)\to \Cl(R)$.  
That this map is surjective follows from a 
standard exercise;  see for instance~\cite[(1.4.18)]{bruns:cmr}.
For the displayed relations it suffices to show, for each $\fa,\fb\in D(R)$
$$R:_Q(R:_Q\fa\fb)\cong (\fa\otimes_R\fb)^{\star\star}
\qquad
\fa:_Q \fb\cong\Hom_R(\fb,\fa).$$ 
(The third condition is obvious from the definition of $f$.) 
The second of these has already been discussed.
For the first, it suffices to show 
$$(\fa\fb)^{\star\star}\cong(\fa\otimes_R\fb)^{\star\star}.$$
Let $K$ be the kernel of the multiplication map $\mu\colon\fa\otimes_R\fb\to\fa\fb$.
One checks that the map $\mu\otimes_R K$ is an isomorphism, and so
$K$ is torsion.  It follows that $K^{\star}=0$ and so
$(\fa\fb)^{\star}\cong(\fa\otimes_R\fb)^{\star}$
and $(\fa\fb)^{\star\star}\cong(\fa\otimes_R\fb)^{\star\star}$.

It follows readily from the isomorphisms described above that, if $\p\in Z(R)$ and $\ell>0$,
then $\ell[\p]=[\p^{(\ell)}]$ in $\Cl(R)$.
\end{para}

\section{Semidualizing modules as divisor classes} \label{sec2}

The following 
proposition compares directly to the ``classical'' result for the dualizing 
module which is the prime motivation for our techniques;  
see, e.g., \cite[(3.3.18)]{bruns:cmr}.  
Recall that a 
finite $R$-module 
$M$ has \emph{rank} (respectively, \emph{rank $r$}) if $M_{\p}$ is 
free (respectively, free of rank $r$) over $R_{\p}$ for each 
$\p\in\ass(R)$.  
Of course, condition~\eqref{item06} is satisfied if $R$ is a domain.
Also, using $C=R$ one sees that the
implication \eqref{item07}$\implies$\eqref{item06}
fails in general; see Proposition~\ref{prop03}\eqref{item10}.

\begin{prop} \label{prop02}
Let $C$ be a semidualizing $R$-module, 
and consider the following conditions.
\begin{enumerate}[\quad\rm(i)]
\item \label{item06}
For each $\p\in\ass(R)$, the localization $R_{\p}$ is Gorenstein.
\item \label{item07}
$C$ has rank 1;
\item \label{item08}
$C$ has rank;
\item \label{item09a}
$C$ is isomorphic to an ideal of $R$;
\item \label{item09}
$C$ is isomorphic to an ideal $\fa$ of $R$ with torsion quotient $R/\fa$.
\end{enumerate}
The implications 
\eqref{item06}$\implies$\eqref{item07}$\iff$\eqref{item08}$\iff$\eqref{item09a}$\iff$\eqref{item09}
hold.
\end{prop}

\begin{proof}
\eqref{item06}$\implies$\eqref{item07}.  
For each associated prime $\p$, 
the ring $R_{\p}$ is Gorenstein and 
therefore the semidualizing $R_{\p}$-module
$C_{\p}$ is isomorphic to $R_{\p}$ by~\cite[(8.6)]{christensen:scatac}.

\eqref{item07}$\implies$\eqref{item08} is trivial.  For the converse, 
since $C_{\p}$ is semidualizing for $R_{\p}$, 
it is routine to check that, if $C_{\p}$ 
is free over $R_{\p}$, then it is free of rank 1.

\eqref{item07}$\iff$\eqref{item09a}$\iff$\eqref{item09}.  It is straightforward to 
show that the semidualizing module $C$ is torsion-free;  in fact, 
$\ass(R)=\ass_R(C)$.  The desired biimplications now follow from a 
standard exercise;  see for instance~\cite[(1.4.18)]{bruns:cmr}.
\end{proof}

A \emph{semidualizing ideal} is an ideal that
is semidualizing as an $R$-module.  One consequence of 
Proposition~\ref{prop02} is that, when $R$ is a domain, every semidualizing 
module is isomorphic to a semidualizing ideal.  The next result
provides basic properties of such ideals;  it
compares directly to~\cite[(3.3.18)]{bruns:cmr}.  We 
restrict our attention to proper ideals as the case $\fa=R$ is tedious.
Since a principal ideal generated by a non-zerodivisor  is semidualizing,
but is dualizing if and only if $R$ is Gorenstein, 
the implication (iii)$\implies$(ii)
fails in general.

\begin{prop} \label{prop03}
Let $R$ be a Cohen-Macaulay ring of dimension $d$ and $\fa$ a proper
semidualizing ideal of $R$.  
\begin{enumerate}[\quad\rm(a)]
\item \label{item05}
$\Ht(\fa)=1$ and $R/\fa$ is Cohen-Macaulay of dimension $d-1$.
\item \label{item14}
The quotient $R/\fa$ has $\text{G}_{\fa}$-dimension 1 and
there are isomorphisms
$$\ext^i_R(R/\fa,\fa)\cong\begin{cases}
R/\fa & \text{if $i=1$} \\ 0 & \text{otherwise.}\end{cases} $$
\item \label{item10}
Consider the following conditions:
\begin{enumerate}[\quad\rm(i)]
\item \label{item11}
The quotient $R/\fa$ is a Gorenstein ring;
\item \label{item12}
The ideal $\fa$ is dualizing for $R$;
\item \label{item13}
$R$ is generically Gorenstein.
\end{enumerate}
The implications (i)$\iff$(ii)$\implies$(iii) hold.
\end{enumerate}
\end{prop}

\begin{proof}
The proof of~\eqref{item05} is nearly identical to that 
of~\cite[(3.3.18.b)]{bruns:cmr}, so we omit it here.
For part~\eqref{item14}, use the exact sequence 
\begin{equation}
0\to\fa\to R\to R/\fa\to 0 \label{eq01} \tag{\ref{prop03}.1}
\end{equation}
with the fact that $\fa$ and $R$ are both 
totally $\fa$-reflexive 
to conclude that $R/\fa$ has $\text{G}_{\fa}$-dimension at most 1.
In particular, $\ext_R^i(R/\fa,\fa)=0$ for $i> 1$.
Furthermore, since $\fa$ has rank, 
it contains an element that is both $R$-regular and $\fa$-regular, and
thus $\Hom_R(R/\fa,\fa)=0$.  
Applying $\Hom_R(-,\fa)$ 
to~\eqref{eq01} supplies the exact sequence
\[ 0\to\underbrace{\Hom_R(R,\fa)}_{\fa}
\to\underbrace{\Hom_R(\fa,\fa)}_R\to\ext_R^1(R/\fa,\fa)\to 0 \]
which yields an isomorphism $\ext_R^1(R/\fa,\fa)\cong R/\fa$ and 
the desired conclusions. 

For part~\eqref{item10} we may assume that $R$ is local.
In the following sequence of formal equalities of series,
the first is by~\cite[(1.6.7)]{christensen:scatac} and the third is 
standard
\[ I_R^{\fa}(t)= I_{R/\fa}^{\rhom_R(R/\fa,\fa)}(t)=
I_{R/\fa}^{\shift^{-1}R/\fa}(t)=
t\cdot I_{R/\fa}^{R/\fa}(t) \]
while the second is 
a consequence of part~\eqref{item14}.
It follows that $\id_R(\fa)$ and $\id_{R/\fa}(R/\fa)$ 
are simultaneously finite.  This gives the equivalence of
(i) and (ii), and the 
implication
(ii)$\implies$(iii) is part 
of~\cite[(3.3.18)]{bruns:cmr}.
\end{proof}

The next result simplifies the computation of
$[\fa]+[\fb]$ in $\Cl(R)$ for certain semidualizing modules
$\fa,\fb$ and is a key tool for the proof of Theorem~\ref{thm02}.

\begin{prop} \label{prop03a}
Let  $\fa$ and $\fb$ be semidualizing ideals 
such that $\fa\otimes_R\fb$ is semidualizing.  
The natural multiplication 
map $\fa\otimes_R\fb\to\fa\fb$ is an isomorphism.
\end{prop}

\begin{proof}
The map $\fa\otimes_R\fb\to\fa\fb$ is always surjective, so it remains 
to verify injectivity.  Let $U$ denote the compliment in $R$ of the 
union of the associated primes of $R$.  Since $\fa$ and $\fb$ have 
rank, the same is true of $\fa\otimes_R\fb$.  Furthermore, the fact 
that $\fa\otimes_R\fb$ is semidualizing implies that 
$\fa\otimes_R\fb$ is torsion-free.  This 
yields the injectivity of the localization map
$\fa\otimes_R\fb\to U^{-1}(\fa\otimes_R\fb)$ in the following commuting diagram
where the maps $(1)$ and $(2)$ are given by the 
appropriate multiplication and the others are the natural ones.
\[ \xymatrix{
\fa\otimes_R\fb \ar[r]^{(1)} \ar@{^{(}->}[d] & \fa\fb \ar@{^{(}->}[r]
& R \ar@{^{(}->}[r] & U^{-1}R \ar[d]^= 
\\
U^{-1}(\fa\otimes_R\fb) \ar[r]^-{\cong} & U^{-1}\fa\otimes_{U^{-1}R}U^{-1}\fb 
\ar@{^{(}->}[r]^-{(2)} & (U^{-1}\fa) (U^{-1}\fb) \ar@{^{(}->}[r]
 & U^{-1}R
} \] 
The map $(2)$ is injective, since $U^{-1}\fa$ and 
$U^{-1}\fb$ are $U^{-1}R$-free of rank 1.  It follows that the map $(1)$ 
must be injective, as desired.
\end{proof}

The next result
supplies the main tool for this investigation.
Note that Theorem~\ref{thm02} shows that $\s_0(R)$ cannot 
be given a group structure making the inclusion 
into a group homomorphism.

\begin{prop} \label{thm01}
Let $R$ be a normal domain.
Each semidualizing $R$-module $C$ is  a 
rank 1 reflexive module,  
so there is a natural inclusion
$\s_0(R)\subseteq\Cl(R)$.
\end{prop}

\begin{proof}
It suffices to verify the first statement.
Proposition~\ref{prop02} shows that $C$ has 
rank 1.  For each prime ideal $\p$ of height 1, the ring $R_{\p}$ is 
regular as $R$ is $(R_1)$, and so $C_{\p}\cong R_{\p}$.  
Since $R$ is $(S_2)$ and 
$\depth_{R_{\p}}(C_{\p})=\depth(R_{\p})$ for each prime ideal $\p$, 
the reflexivity of $C$ follows 
from~\cite[(1.4.1)]{bruns:cmr}.  
\end{proof}

We record an immediate corollary.

\begin{cor} \label{cor01}
Every normal domain with finite divisor class group is 
$\s_0$-finite, and
every Cohen-Macaulay normal domain with finite divisor class group is 
$\s$-finite.\qed
\end{cor}

Since a Cohen-Macaulay normal domain $R$ with $\Cl(R)=0$ is 
Gorenstein, 
we note that there are non-Gorenstein rings 
that satisfy the hypotheses of the corollary.  For instance, if 
$k$ is a field and $\mathbf{X}$  
a symmetric $n\times n$ matrix of variables and $r$ an integer 
such that $0<r<n$, then the ring 
$R=k[\mathbf{X}]/I_{r+1}(\mathbf{X})$ is a Cohen-Macaulay normal domain
with $\Cl(R)\cong \zz/(2)$ and is non-Gorenstein
if and only if $r\equiv n \pmod 2$.  
Here $I_{r+1}(X)$ is the ideal generated 
by the minors of $\mathbf{X}$ of size $r+1$;  see~\cite[(7.3.7.c)]{bruns:cmr}.
Determinantal rings will be of particular interest in Section~\ref{sec3}.

Proposition~\ref{thm01} points toward a plethora of examples of nonlocal
rings that are neither $\s$-finite nor $\s_0$-finite.
Hence the local hypothesis in Question~\ref{q01}.

\begin{para} \label{para02}
the \emph{Picard group} of a normal domain $R$, denoted $\Pic(R)$,
is the set of
isomorphism classes of finitely generated locally free (i.e.,
projective) $R$-modules of rank 1 with operation given by tensor product. 
The inverse of an element $[P]\in\Pic(R)$ is the class $[P]^{-1}=[\Hom_R(P,R)]$;
see~\cite[p.~105]{fossum:dcgkd}.
It is straightforward to show that
there are natural 
inclusions $\Pic(R)\subseteq\s_0(R)\subseteq\Cl(R)$ for any normal 
domain $R$.
Using~\cite[(3.2)]{frankild:appx} one sees that the first inclusion is an
equality when $R$ is Gorenstein.
Each inclusion is an equality when $R$ is a Dedekind 
domain by Fossum~\cite[(18.5)]{fossum:dcgkd}.  

A result of 
Claborn~\cite[(14.10)]{fossum:dcgkd} states that any Abelian 
group $G$ can be realized as the divisor class group of a Dedekind 
domain.  
In particular, for any 
Abelian group $G$, regardless of the cardinality, there is a Dedekind 
domain $R$ such that the sets $\s(R)$ and $\s_0(R)$ are in bijection with $G$.
\end{para}

We observe that~\cite[(3.1.b)]{frankild:appx}
implies that the hypothesis of the next result is
satisfied when $C$ is 
$C''$-reflexive and $C'=\Hom_R(C,C'')$.  
Compare to Proposition~\ref{prop03a}.

\begin{prop} \label{prop03b}
Let $R$ be a normal domain and $C,C'$ semidualizing 
$R$-modules.
If $C'\otimes_R C$ is $R$-semidualizing,  then 
$[C'\otimes_R C]=[C']+[C]$ in $\Cl(R)$.
\end{prop}

\begin{proof}
Since
$C'\otimes_R C$ is semidualizing, it is reflexive, so 
$(C'\otimes_R C)^{\star\star}\cong C'\otimes_R C$, and 
so $[C'\otimes_R C]=(C'\otimes_R C)^{\star\star}=[C']+[C]$;
see~\ref{paranew7}.
\end{proof}

The inclusion $\s_0(R)\subseteq\Cl(R)$
is well-behaved with respect to certain operations that 
are defined on both sets.  The remainder of this section is devoted to
describing some of this 
behavior.  We begin by describing base-change maps for
Picard groups and sets of semidualizing objects.

\begin{para} \label{para03}
Let $\vf\colon R\to S$ be a ring homomorphism. 
\begin{enumerate}[\quad(a)]
\item \label{item92d}
The assignment 
$L\mapsto L\lotimes_R S$
yields a well-defined group 
homomorphism $\Pic(\vf)\colon\Pic(R)\to\Pic(S)$; see~\cite[discussion after (18.3)]{fossum:dcgkd}.  
\item \label{item92a}
When $\fd(\vf)$ is finite
and $K$ is a semidualizing complex, the $S$-complex $K\lotimes_R S$ 
is semidualizing by~\cite[(4.5)]{frankild:appx}, 
and the 
assignment $K\mapsto K\lotimes_R S$ gives rise 
to a well-defined order-respecting
map $\s(\vf)\colon \s(R)\to\s(S)$ by~\cite[(4.7)]{frankild:appx}.
\item \label{item92b}
When $\fd(\vf)$ is finite
and $C$ is 
a semidualizing $R$-module, the $S$-module $C\otimes_R S$ is 
semidualizing and $\tor^R_{\geqslant 1}(C,S)=0$ by~\cite[(4.5)]{frankild:appx};
the assignment $C\mapsto C\otimes_R S$ induces
a well-defined order-respecting
map  $\s_0(\vf)\colon \s_0(R)\to\s_0(S)$
by~\cite[(4.7)]{frankild:appx}. 
\end{enumerate}
\end{para}

Next we consider the divisor class group.
Part~\eqref{item92c} in
the following lemma is well-known, but we include it here for completeness;
see~\cite[Sec.~6]{fossum:dcgkd}.

\begin{lem} \label{para03x}
Let $\vf\colon R\to S$ be a homomorphism of finite flat dimension between normal
domains.  
\begin{enumerate}[\quad\rm(a)]
\item \label{item9003y}
Fix $\q\in\spec(S)$  and set $\p=\vf^{-1}(\q)$.  If $\Ht(\q)\leq 1$, then $R_{\p}$ is regular.
\item \label{item9003x}
If $M$ is a reflexive $R$-module of rank 1, then $\rank_S(M\otimes_R S)=1$.
\item \label{item90}
If 
$\vf$ is module-finite,
the assignment 
$M\mapsto \Hom_S(\Hom_S(M\lotimes_R S,S),S)$
yields a well-defined group 
homomorphism $\Cl(\vf)\colon\Cl(R)\to\Cl(S)$.  
\item \label{item92c}
If 
$\vf$ is flat, then 
the assignment 
$M\mapsto M\lotimes_R S$
yields a well-defined group 
homomorphism $\Cl(\vf)\colon\Cl(R)\to\Cl(S)$.  
\end{enumerate}
\end{lem}

\begin{proof}
\eqref{item9003y}
The induced map
$\vf_{\q}\colon R_{\p}\to S_{\q}$ has finite flat dimension.  Since $S$ 
is normal, it satisfies Serre's 
condition ($\text{R}_1$) and so the local ring $S_{\q}$ is regular.  
It follows from~\cite[Thm.~R]{apassov:afm} that $R_{\p}$ is also regular.

\eqref{item9003x}
To show that $M\otimes_R S$ has rank 1, it suffices to set $\q=(0)S$ and exhibit an isomorphism
$(M\otimes_R S)_{\q}\cong S_{\q}$.  
With $\p=\Ker(\vf)$, part~\eqref{item9003y}
implies $\Cl(R_{\p})=0$ and so $M_{\p}\cong R_{\p}$.
The next isomorphisms now follow readily
$$
(M\otimes_R S)_{\q}\cong M_{\p}\otimes_{R_{\p}}S_{\q}\cong R_{\p}\otimes_{R_{\p}}S_{\q}\cong S_{\q}
$$
and provide the desired conclusion.

\eqref{item90}
Using part~\eqref{item9003y}, this follows from~\cite[(1.2.1)]{spiroff:diss}.

\eqref{item92c}
For a finitely generated $R$-module $U$, one has a natural $S$-linear map
\begin{align*}
f_{U}^{}\colon \Hom_R(U,R)\otimes_R S&\to\Hom_S(U\otimes_RS,S) \\
\psi\otimes s&\mapsto[u\otimes s'\mapsto\vf(\psi(u))ss']
\end{align*}
which is readily seen to be an isomorphism because $\vf$ is flat.

Let $M$ and $N$ be rank 1 reflexive $R$-modules.  
The $S$-module
$M\otimes_R S$ has rank 1 by part~\eqref{item9003x}.  The flatness of $\vf$
provides the following commutative diagram
from which one concludes that $M\otimes_R S$ is a reflexive $S$-module.
$$
\xymatrix{
M\otimes_R S \ar[rr]_-{\cong}^-{b^R_M\otimes_R S} \ar[d]_{b^S_{M\otimes_R S}}
&& \Hom_R(\Hom_R(M,R),R)\otimes_R S \ar[d]^{\cong}_{f^{}_{\Hom_R(M,R)}} \\
\Hom_S(\Hom_S(M\otimes_R S,S),S) \ar[rr]^{\Hom_S(f_M^{},S)}_{\cong}
&& \Hom_S(\Hom_R(M,R)\otimes_R S,S)
}
$$
Hence, the map $\Cl(\vf)$ is well-defined.
The fact that $\Cl(\vf)$ is a group homomorphism
follows from the next sequence of isomorphisms 
\begin{align*}
\Hom_R(\Hom_R(M\otimes_R N,R),R)\otimes_R S 
& \stackrel{(1)}{\cong} \Hom_S(\Hom_R(M\otimes_R N,R)\otimes_R S,S) \\
& \hspace{-.75cm}\stackrel{(2)}{\cong} \Hom_S(\Hom_S((M\otimes_R N)\otimes_R S,S),S) \\
& \hspace{-1.5cm}\stackrel{(3)}{\cong} \Hom_S(\Hom_S((M\otimes_R S)\otimes_S(N\otimes_R S),S),S)
\end{align*}
where (1) is $f^{}_{\Hom_R(M\otimes_R N,R)}$, (2) is $\Hom_S(f^{}_{M\otimes_R N},S)$,
and (3) is standard.
\end{proof}

\begin{lem} \label{prop04}
Let $\vf\colon R\to S$ be a homomorphism of finite flat dimension.
\begin{enumerate}[\quad\rm(a)]
\item \label{item15}
If $R$ and $S$ are normal domains and $\vf$ is module-finite or flat,
then the following diagram commutes.
\[ \xymatrix{
\s_0(R) \ar@{^{(}->}[r] \ar[d]_{\s_0(\vf)} & \Cl(R) \ar[d]_{\Cl(\vf)} \\
\s_0(S) \ar@{^{(}->}[r] & \Cl(S) 
} \]
In particular, if $\Cl(\vf)$ is injective, then so is $\s_0(\vf)$.
\item \label{item16}
Assume that the image of the map $\spec(\vf)\colon
\spec(S)\to\spec(R)$ contains $\mspec(R)$.
If the map $\Pic(\vf)\colon\Pic(R)\to\Pic(S)$
is injective, e.g., if $\vf$ is surjective or local,
then 
$\s_0(\vf)$ and $\s(\vf)$ are also injective.  
\item \label{item19}
Assume that $R,S$ are normal domains and $\vf$ is faithfully flat.
If $\Cl(\vf)$ is surjective, then so is
$\s_0(\vf)$.
\item \label{item19a}
Assume that $R,S$ are normal domains and $\vf$ is faithfully flat.
If $\Cl(\vf)$ is bijective, then 
$\s_0(\vf)$ is a perfectly order-respecting bijection. 
\end{enumerate}
\end{lem}

\begin{proof}
When $\vf$ is flat, the commutativity of the diagram in~\eqref{item15} follows readily from
the definitions.
When $\vf$ is module finite and $C$ is a semidualizing $R$-module,
the fact that
$C\otimes_R S$ is semidualizing for $S$ implies that it is reflexive, and the
commutativity of the diagram follows easily.
Part~\eqref{item16}
is contained in~\cite[(4.9),(4.11)]{frankild:appx},
and~\eqref{item19}
is in~\cite[(4.5)]{frankild:appx}.
When $\Cl(\vf)$ is bijective, the map $\Pic(\vf)$ is injective, 
so~\eqref{item19a} follows from parts~\eqref{item16} and~\eqref{item19}
with~\cite[(4.8)]{frankild:appx}.
\end{proof}

When $\vf\colon R\to S$
is a local homomorphism of finite flat dimension, 
it is a straightforward exercise to show that, if $\s(\vf)$ is
a perfectly order-respecting bijection, then it is also an isometry
with respect to the metric structure defined in~\cite{frankild:sdcms}.
For instance, this holds under the hypotheses of Lemma~\ref{prop04}\eqref{item19a}
when $S$ is Cohen-Macaulay.
Some particular instances of this are provided in the next corollary.
Others are given in  Corollary~\ref{cor03}
and Proposition~\ref{prop12}.

\begin{cor} \label{cor02}
Let $R$ be a normal domain and $\mathbf{X}=X_1,\ldots,X_n$ 
a sequence of
variables.  For the following flat $R$-algebras $S$, the map
$\s_0(R)\to\s_0(S)$ is a perfectly order-respecting bijection:
\begin{enumerate}[\quad\rm(a)]
\item \label{item17}
$S=R[\mathbf{X}]$;
\item \label{item17a}
$S=R[\mathbf{X}][f_1^{-1},\ldots,f_i^{-1}]$ where
$f_1,\ldots,f_i$ are prime elements of $R[\mathbf{X}]$ and the 
ring homomorphism
$R\to R[\mathbf{X}][f_1^{-1},\ldots,f_i^{-1}]$ is faithfully flat;
\item \label{item17c}
$S=R[\mathbf{X}]_{\m R[\mathbf{X}]}$ when
$R$ is local with maximal 
ideal $\m$;
\item \label{item18a}
$S=R[\![\mathbf{X}]\!]_{\m R[\![\mathbf{X}]\!]}$ when
$R$ is local with maximal ideal $\m$ and the $\m$-adic completion 
of $R$ is normal.
\end{enumerate}
\end{cor}

\begin{proof}
By Lemma~\ref{prop04}\eqref{item19a}, it suffices to note that the maps
$\Cl(R)\to\Cl(S)$
are bijective; 
see~\cite[(7.3),(8.1),(8.9),(19.15)]{fossum:dcgkd}.
\end{proof}

The following is an important case when localization induces a 
bijection on the set of semidualizing modules.

\begin{prop} \label{prop05}
Let $R=\coprod_{i\geq 0}R_i$ be a graded normal domain 
such that $(R_0,\m_0)$ is local. 
Setting $\m=\m_0+\coprod_{i\geq 1}R_i$,
the natural 
map 
$\s_0(R)\to\s_0(R_{\m})$ is 
a perfectly order-respecting bijection.
\end{prop}

\begin{proof}
Let $\vf\colon R\to R_{\m}$ be the localization map.
Using~\cite[(2.14)]{frankild:appx},
the argument
of~\cite[(10.3)]{fossum:dcgkd} shows that $\Cl(\vf)$
is bijective.
From Lemma~\ref{prop04}\eqref{item15} it follows that 
$\s_0(\vf)$ is injective.  To show surjectivity, fix a semidualizing 
$R_{\m}$-module $L$.  Use the surjectivity of $\Cl(\vf)$ 
and~\cite[(10.2)]{fossum:dcgkd} to obtain a homogeneous reflexive 
ideal $\fa$ of $R$ such that $\fa_{\m}\cong L$.  Since $L$ is 
$R_{\m}$-semidualizing, the $R$-module $\fa$ is $R$-semidualizing
by~\cite[(2.15.a)]{frankild:appx},
and it follows that
$\s_0(\vf)$ is bijective.  The fact that $\s_0(\vf)$ is
perfectly order-respecting then follows from~\cite[(2.15.b)]{frankild:appx}.
\end{proof}

If $R$ is a local ring with completion map $\vf\colon R\to
\comp{R}$, then the map $\s_0(\vf)$ is not usually surjective.  
Indeed, there exist a Cohen-Macaulay local 
ring $R$ that does not admit a dualizing 
module; the complete local ring $\comp{R}$ does admit a dualizing 
module $\omega$, and it is straightforward to show that 
$[\omega]\in\s_0(\comp{R})$ cannot be in the image of $\s_0(\vf)$.
However, a result of Flenner~\cite[(1.4)]{flenner:1} 
can be applied in certain cases to provide bijectivity.
See 
Corollary~\ref{cor08}\eqref{item72} for a generalization,
and also~\cite[(1.1)]{christensen:dscrap}.
Recall that a ring is \emph{super-normal} if it satisfies Serre's 
conditions $(S_3)$ and $(R_2)$.
Further, note that a ring $R$ satisfying the hypotheses of the next result is excellent
because it is finitely generated over the complete local ring $R_0$.
In particular, the complete ring $\comp{R}$ is also a super-normal domain.

\begin{cor} \label{cor03}
Let $R=\coprod_{i\geq 0}R_i$ be a graded super-normal domain with 
$(R_0,\m_0)$ local and complete, and set $\m=\m_0+\coprod_{i\geq 1}R_i$ and
$\comp{R}=\prod_{i\geq 0}R_i$.
The induced maps $\s_0(R)\to\s_0(\comp{R})$ 
and $\s_0(R_{\m})\to\s_0(\comp{R})$ are 
perfectly order-respecting bijections.
\end{cor}

\begin{proof}
The ring $\comp{R}$ is the $\m$-adic completion of 
$R_{\m}$, and since $R$ is excellent and super-normal, the same is 
true of $R_{\m}$ and $\comp{R}$.
Let $\vf\colon R\to R_{\m}$ be the localization map and $\psi\colon
R_{\m}\to \comp{R}$ the completion map.   
By Proposition~\ref{prop05}, the map $\s_0(\vf)$ is 
a perfectly order-respecting bijection, 
so the equality
$\s_0(\psi\vf)=\s_0(\psi)\s_0(\vf)$ 
shows that we need only verify the same for $\s_0(\psi)$.  Since
$\psi$ is flat and local, Lemma~\ref{prop04}\eqref{item16} 
supplies the
injectivity of $\s_0(\psi)$.  The surjectivity 
is  a consequence of Lemma~\ref{prop04}~\eqref{item19}, 
as~\cite[(1.4)]{flenner:1} guarantees that 
$\Cl(\psi\vf)=\Cl(\psi)\Cl(\vf)$ 
is surjective, and therefore that $\Cl(\psi)$ is surjective.  
\end{proof}

Here is 
the first indication that our methods have applications outside the 
normal domain arena.
See Corollary~\ref{cor09} for a more general statement.

\begin{cor} \label{cor08}
With $R,\m,\comp{R}$ as in Corollary~\ref{cor03}, fix
an $R_{\m}$-regular sequence $\mathbf{y}=y_1,\ldots,y_q\in\m R_{\m}$. 
\begin{enumerate}[\quad\rm(a)]
\item \label{item72}
The natural homomorphisms
$R_{\m}\to R_{\m}/(\mathbf{y})\to\comp{R}/(\mathbf{y})$
induce 
perfectly order-respecting bijections
\[ 
\s_0(R_{\m})\xra{\approx}
\s_0(R_{\m}/(\mathbf{y}))\xra{\approx}
\s_0(\comp{R}/(\mathbf{y})). \]
\item \label{item72a}
Let $R'$ denote either
$R_{\m}/(\mathbf{y})$ or $\comp{R}/(\mathbf{y})$ and fix an $\m R'$-primary
ideal $J\subset R'$.  
If $C'$ is a semidualizing $R'$-module, then $e(J,C)=e(J,R')$.
\item \label{item73}
If $\mathbf{y}$ is an $R$-regular sequence in $\m$, then the composition
of induced maps
\[ \s_0(R)\to\s_0(R/(\mathbf{y}))\to\s_0(R_{\m}/(\mathbf{y})) \]
is a perfectly order-respecting bijection;  
thus, the first map is a perfectly order-respecting injection and the 
second is surjective.
\end{enumerate}
\end{cor}

\begin{proof}
\eqref{item72}
The rings under consideration fit into the commutative diagram 
of local ring homomorphisms on the left
$$ \xymatrix{
R_{\m} \ar[r]^{\beta} \ar[d]_{\alpha_{\m}} & \comp{R} \ar[d]_{\comp{\alpha}} 
& \s_0(R_{\m}) \ar[r]_{\approx}^{\s_0(\beta)} \ar[d]_{\s_0(\alpha_{\m})} 
& \s_0(\comp{R}) \ar[d]^{\approx}_{\s_0(\comp{\alpha})}
\\
R_{\m}/(\mathbf{y}) \ar[r]^{\beta'} & \comp{R}/(\mathbf{y})
& \s_0(R_{\m}/(\mathbf{y})) \ar[r]^{\s_0(\beta')} & \s_0(\comp{R}/(\mathbf{y}))
} $$
and the second commutative diagram arises by applying
$\s_0(-)$ to the first.
The maps $\s_0(\beta)$ and $\s_0(\comp{\alpha})$ are perfectly 
order-respecting bijections
by 
Corollary~\ref{cor03} and~\cite[(4.2)]{frankild:sdcms}, respectively;  
see also Gerko~\cite[(3)]{gerko:osmagi}. 
From the diagram, it follows that $\s_0(\beta')$ is surjective, and so is 
bijective by Lemma~\ref{prop04}\eqref{item16}.  
That it is perfectly order-respecting 
is then a consequence of~\cite[(4.8)]{frankild:appx}.
Thus, $\s_0(\alpha_{\m})$ is a  perfectly order-respecting bijection, as well.

\eqref{item72a}
Set $\wt{R}=R_{\m}$ or $\comp{R}$, according to whether
$R'=R_{\m}/(\mathbf{y})$ or $\comp{R}/(\mathbf{y})$.
Fix a semidualizing $\wt{R}$-module
$C$ such that $C'\cong C\otimes_{\wt{R}}R'$.
Since $\wt{R}$ is an excellent local domain, 
one applies Lemma~\ref{lem:mult}\eqref{item:mult:b}
to obtain the desired conclusion.

\eqref{item73}
When $\mathbf{y}$ is an $R$-regular sequence in $\m$, 
there is a commutative diagram
\[ \xymatrix{ 
\s_0(R) \ar[r]_{\approx}^{\s_0(\gamma)} \ar[d]_{\s_0(\alpha)} 
& \s_0(R_{\m}) \ar[d]^{\approx}_{\s_0(\alpha_{\m})}  \\
\s_0(R/(\mathbf{y})) \ar[r]^{\s_0(\gamma')} & \s_0(R_{\m}/(\mathbf{y})) 
} \]
where $\s_0(\gamma)$ and $\s_0(\alpha_{\m})$ are perfectly
order-respecting bijections by Proposition~\ref{prop05}
and part~\eqref{item72}.
The injectivity of $\s_0(\alpha)$ and surjectivity of 
$\s_0(\gamma')$ follow immediately.
To see that $\s_0(\alpha)$ 
is perfectly order-respecting, fix $[C],[C']\in\s_0(R)$
such that $\s_0(\alpha)([C])\tri\s_0(\alpha)([C'])$.
It follows that
\begin{align*}
\s_0(\alpha_{\m})(\s_0(\gamma)([C]))
&=\s_0(\gamma)(\s_0(\alpha)([C]))\\
&\tri\s_0(\gamma)(\s_0(\alpha)([C']))\\
&=\s_0(\alpha_{\m})(\s_0(\gamma)([C']))
\end{align*}
and so $[C]\tri[C']$ since
$\s_0(\alpha_{\m})$ and $\s_0(\gamma)$ are perfectly order-respecting.
\end{proof}

The surjectivity of the natural map 
$\s_0(R_{\m})\to \s_0(R_{\m}/(\mathbf{y}))$ 
in the 
corollary does not hold 
for more general local rings.  Indeed, let $(A,\n)$ be a local 
Cohen-Macaulay ring that does not admit a dualizing module.  
If $\mathbf{y}\in\n$ is a system of parameters of $A$, then the map
$\s_0(A)\to\s_0(A/(\mathbf{y}))$ is not surjective, as 
$A/(\mathbf{y})$ is Artinian
and therefore admits a dualizing module. See~\cite[(5.5)]{christensen:dscrap}
for further discussion.

When the homomorphism $\vf\colon R\to S$ is part of a retract pair,
the next result sometimes allows one to conclude that $\s_0(\vf)$ is
bijective.  Examples of such retract pairs can be found in power
series and localized polynomial extensions:
\begin{enumerate}[\quad\rm(a)]
\item \label{item74}
The natural maps $R\to R[\![\mathbf{X}]\!]$ and 
$R[\![\mathbf{X}]\!]\to R$.
\item \label{item75}
The
natural maps $R\to R[\mathbf{X}]_{\n}$ and 
$R[\mathbf{X}]_{\n}\to R$, when 
$(R,\m)$ is local 
and $\n=(\m,X_1-a_1,\ldots,X_n-a_n)R[\mathbf{X}]$
for a sequence $a_1,\ldots,a_n\in R$.  
\end{enumerate}
Note that 
the 
rings involved are not assumed to be normal domains, so one cannot 
use the divisor class group directly.  However, the method of proof is 
taken directly from the corresponding results for 
divisor class groups.  

\begin{prop} \label{prop12}
Let $\vf\colon R\to S$  and $\psi\colon S\to R$ be homomorphisms
of finite flat dimension such that the composition $\psi\vf$ is
the identity on $R$.  
If $\ker(\psi)$ is contained in the Jacobson radical of $S$,
then the induced maps $\s_0(\vf),\s_0(\psi),\s(\vf),\s(\psi)$
are perfectly order-respecting bijections.
\end{prop}

\begin{proof}
\eqref{item74} 
Since the composition $\psi\vf\colon R\to R$ is the identity, the 
same is true of the composition $\s_0(\psi)\s_0(\vf)$.  
In particular, $\s_0(\psi)$ is surjective.  Since $\ker(\psi)$ is in the 
Jacobson radical of $S$, Lemma~\ref{prop04}\eqref{item16}
guarantees that this 
map 
is bijective, 
and therefore so is 
$\s_0(\vf)$.  The same argument works for $\s(\psi)$ and $\s(\vf)$.
\end{proof}

This is a surprising departure from the parallels we have seen
between the behavior of $\Cl(-)$ and $\s_0(-)$, as it is known that,
when $R$ is a normal domain, the map 
$\Cl(R)\to\Cl(R[\![\mathbf{X}]\!])$ need
not be bijective;  see Danilov~\cite{danilov:2,danilov:1,danilov:3}.

The proof of Proposition~\ref{prop12}
can be translated easily to show that the natural maps
$R\to R[\mathbf{X}]\to R$ induce injections and surjections respectively
on sets of semidualizing objects.  However, 
we can only say more about these maps when we assume
that $R$ is a normal domain, by using Corollary~\ref{cor02}.

\begin{question} \label{q04}
Must the induced maps $\s_0(R)\to\s_0(R[X])\to\s_0(R)$ and
$\s(R)\to\s(R[X])\to\s(R)$ all be bijective?
\end{question}

\section{Analysis of special cases} \label{sec3}

We begin with some
notation and facts 
on determinantal rings from~\cite{bruns:dr}.

\begin{para}
\label{para:det1}
Let $A$ be a
Noetherian ring and $m,n,r$ nonnegative integers satisfying
$r<\min\{m,n\}$.  If $\mathbf{X}=\{X_{ij}\}$ is an $m\times n$ matrix of
variables, then set 
\[ R=R_{r+1}(A;m,n)=A[\mathbf{X}]/I_{r+1}(\mathbf{X}) \]
where $I_{r+1}(\mathbf{X})$
denotes the ideal generated by the minors of $\mathbf{X}$ of size $r+1$.  
If $A$ is a normal domain
(respectively, is Cohen-Macaulay or 
is $(S_3)$ or is $(R_2)$)
then so is $R$; see~\cite[(6.3),(5.17),(5.16),(6.12)]{bruns:dr}.
The ring $R$ is Gorenstein if and only if 
$A$ is Gorenstein and either $m= n$ or $r= 0$ by~\cite[(8.9)]{bruns:dr}.

Assume that $A$ is a normal domain and $r>0$.
Let $\p$ be the ideal of $R$ generated by the $r$-minors of the first 
$r$ rows of the residue matrix $\mathbf{x}$.  
The ideal $\p$ is prime, and there is an isomorphism $\Cl(R)\cong\Cl(B)\oplus \mathbb{Z}$
where the summand $\mathbb{Z}$ is generated by  $[\p]$;
see~\cite[(8.4)]{bruns:dr}.
For $\ell>0$ one 
has $-[\p]=[\Hom_R(\p,R)]=[\q]$
where $\q$ is the prime ideal of $R$ generated by the 
$r$-minors of the first $r$ columns of $\x$.
For $\ell\ge0$ one has $\ell[\p]=[\p^{(\ell)}]=[\p^{\ell}]$
and $-\ell[\p]=[\q^{(\ell)}]=[\q^{\ell}]$
by~\cite[(7.10)]{bruns:dr}, and we write
$\p^{-\ell}$ in place of $\q^\ell$.
If $A$ is also Gorenstein local and $m\geq n$, then
$R$ admits a unique (up to isomorphism) dualizing module $\omega\cong\p^{(m-n)}=\p^{m-n}$;
see~\cite[(7.10),(8.8)]{bruns:dr}.

Assume that $A$ is a field and $m\geq n\geq r>1$.  Let 
$\mathbf{Y}=\{Y_{pq}\}$ be an $(m-1)\times(n-1)$ matrix of 
variables and set $R'=R_r(A;m-1,n-1)$.  
Let $\p'$ be the ideal of 
$R'$ generated by the $(r-1)$-minors of the first 
$r-1$ rows of the residue matrix $\mathbf{y}$.  
The discussion above implies that
$R'$ has a unique dualizing module, namely, the ideal
$(\p')^{(m-1)-(n-1)}=(\p')^{m-n}$.
We consider three homomorphisms of normal domains
\[ R\xra{\vf} R_{x_{11}}\xleftarrow[\cong]{\rho}
R'[X_{11},\ldots,X_{m1},X_{12},\ldots,X_{1n}]_{X_{11}}
\xleftarrow{\psi} R' \]
where 
$\rho$ is given by $y_{ij}\mapsto x_{i+1 j+1}-x_{1 j+1}x_{i+1,1}x_{11}^{-1}$
and $\vf$ and $\psi$ are the natural flat maps.
By~\cite[(7.3.3)]{bruns:cmr} the map $\rho$ is an isomorphism. 
Further, the
induced maps
are all isomorphisms between groups isomorphic to $\mathbb{Z}$
\[ \Cl(R)\xra{\cong} \Cl(R_{x_{11}})\xleftarrow{\cong} 
\Cl(R'[X_{11},\ldots,X_{m1},X_{12},\ldots,X_{1n}]_{X_{11}})
\xleftarrow{\cong} \Cl(R') \]
and $\Cl(\vf)([\p])=\Cl(\rho\psi)([\p'])$; see Lemma~\ref{para03x}\eqref{item92c} 
and~\cite[proof of (7.3.6)]{bruns:cmr}.
For each integer $\ell\geq 0$, the additivity of $\Cl(\vf)$ and $\Cl(\rho\psi)$ 
provides the second equality in the next sequence
while the others are from the previous paragraph.
$$\Cl(\vf)([\p^{\ell}])=\Cl(\vf)(\ell[\p])=\Cl(\rho\psi)(\ell[\p'])=\Cl(\rho\psi)([(\p')^{\ell}])$$
Lemma~\ref{prop04}\eqref{item19a} implies that $\s_0(\psi)$ 
is bijective.
Let $C$ be a semidualizing $R$-module  and let $c$ be the unique 
integer with $[C]=c[\p]=[\p^c]$ in $\Cl(R)$.  
The $R_{x_{11}}$-module 
$C\otimes_R R_{x_{11}}$ is semidualizing by~\ref{para03}\eqref{item92b}, and
we compute its class in 
$\s_0(R_{x_{11}})$ in the next sequence 
where the equalities follow from 
Lemma~\ref{prop04}\eqref{item15}, the choice of $c$, and
the previous displayed sequence.
$$\s_0(\vf)([C])=\Cl(\vf)([C])=\Cl(\vf)([\p^c])=\Cl(\rho\psi)([(\p')^c])$$
As $\rho\psi$ is flat, this provides an isomorphism of $R_{x_{11}}$-modules
$$(\p')^c\otimes_{R'}R_{x_{11}}\cong C\otimes_RR_{x_{11}}.$$ 
Since $C$ is $R$-semidualizing,
the module $C\otimes_RR_{x_{11}}$ is $R_{x_{11}}$-semidualizing, and
the last isomorphism implies that $(\p')^c$ is $R'$-semidualizing
by~\cite[(4.5)]{frankild:appx}.
\end{para}

\begin{thm} \label{thm02}
Let $k$ be a field and $m,n,r$ nonnegative integers such that
$r<\min\{m,n\}$.  The ring $R=R_{r+1}(k;m,n)$ satisfies
$\s_0(R)=\{[R],[\omega]\}$ 
where $\omega$ is a dualizing module for $R$.  
In particular the cardinality of $\s_0(R)$ is
\[ \card\s_0(R)=\begin{cases} 1 & \text{when $m=n$ or $r=0$} \\
        2 & \text{when $m\neq n$ and $r\neq 0$.} \end{cases} \]
\end{thm}

\begin{proof}
If $r=0$ or $m=n$, then $R$ is Gorenstein and
the result follows from~\ref{para02}.  Assume for the remainder of the proof 
that $r>0$ and $m\neq n$.  
We may also assume that $n\leq m$, as replacing $\mathbf{X}$ with its 
transpose yields an isomorphism $R_{r+1}(k;m,n)\cong 
R_{r+1}(k;n,m)$.  Let $x_{ij}$ denote the residue of $X_{ij}$ in $R$.

We have the containment $\s_0(R)\supseteq\{[R],[\p^{m-n}]\}$ from~\ref{para:det1},
so it remains to verify the containment $\s_0(R)\subseteq\{[R],[\p^{m-n}]\}$.
Let $C$ be a semidualizing $R$-module  and let $c$ be the unique 
integer with $[C]=c[\p]$ in $\Cl(R)\cong\mathbb{Z}[\p]$.  

Following the proof of~\cite[(7.3.6)]{bruns:cmr}, we 
use induction on $r$ to reduce to the case $r=1$.  
Suppose that $r>1$ and employ the notation of the last paragraph of~\ref{para:det1}.
The final conclusion of~\ref{para:det1} says that $(\p')^c$ is $R'$-semidualizing,
and the induction hypothesis states that $\s_0(R')=\{[R'],[(\p')^{m-n}]\}$.
Hence, either $c=0$ or $c=m-n$.  By our choice of $c$, this implies
either $[C]=[\p^0]=[R]$ or $[C]=[\p^{m-n}]$, as desired.

Assume $r=1$.  
As above, it suffices to
show that $c=0$ or $c=m-n$.  
The ring $R$ is a standard graded 
ring over a field and $\p$ is a homogeneous prime ideal. 
For each $v\ge 0$ the power $\p^v$ 
is homogeneous, and so we may 
speak of its minimal number of generators, 
denoted 
$\beta_0(\p^{v})$.
As is noted in~\cite[(9.20)]{bruns:dr}, 
the homogeneous minimal generators of 
$\p^{v}$ are in bijection with the monomials of degree $v$ in the 
ring $k[Z_1,\ldots,Z_n]$.
Since $n\geq 2$,
a routine argument shows
\begin{equation}\label{eq07}\tag{\ref{thm02}.1}
\beta_0(\p^{u})\beta_0(\p^{v})>\beta_0(\p^{u+v})>\beta_0(\p^{v}).
\end{equation}
Also, when $u\leq v$, 
the following sequence implies $\p^{v-u}\cong\Hom_R(\p^{u},\p^{v})$:
$$[\p^{v-u}]=(v-u)[\p]=v[\p]-u[\p]=[\p^v]-[\p^u]=[\Hom_R(\p^u,\p^v)].$$

Suppose that $0<c<m-n$.  
Then $C\cong\p^{c}$,
and
$\p^{m-n}$ is a dualizing module for $R$.  It follows 
from~\cite[(2.14.a),(3.1.a)]{frankild:appx} 
and~\cite[(3.4.a)]{christensen:scatac}
that
$\Hom_R(\p^{c},\p^{m-n})\cong\p^{m-n-c}$ is semidualizing.  
Furthermore, Proposition~\ref{prop03a}  yields
an isomorphism
\[ \p^{c}\otimes_R\p^{m-n-c}\xra{\cong} 
\p^{m-n} \]
and thus the equality
$\beta_0(\p^{m-n})=\beta_0(\p^{c})\beta_0(\p^{m-n-c})$,
contradicting~\eqref{eq07}.  

Next, suppose that $c>m-n$.  As above, we have
\[ \beta_0(\p^{c})>\beta_0(\p^{m-n})
=\beta_0(\p^{c})\beta_0(\Hom_R(\p^{c},\p^{m-n}))
>\beta_0(\p^{c}) \] 
again yielding a 
contradiction.

Finally, suppose that $c<0$.  Then 
$\Hom_R(C,\p^{m-n})\cong\p^{m-n-c}$ is semidualizing.  However, $c<0$ 
implies that $m-n-c>m-n$, contradicting
the previous case.  
\end{proof}

Next, we present the local analogue of Theorem~\ref{thm02}.

\begin{cor} \label{cor04}
With notation as in Theorem~\ref{thm02}, let $\m$ denote the maximal 
ideal of $R$ generated by the residues of the variables $X_{ij}$.
If $R'$ is either the localization $R_{\m}$ or its $\m$-adic 
completion $\comp{R}$, then
$\s_0(R')=\{[R'],[\omega']\}$ 
where $\omega'$ is a dualizing module for $R'$.  
In particular the cardinality of $\s_0(R')$ is
\[ \card \s_0(R')=\begin{cases} 1 & \text{when $m=n$ or $r=0$} \\
        2 & \text{when $m\neq n$ and $r\neq 0$.} \end{cases} \]
\end{cor}

\begin{proof}
The ring $R$ satisfies the hypotheses of Corollary~\ref{cor03} as it 
is Cohen-Macaulay and $(R_2)$ by~\cite[(6.12)]{bruns:dr}.
\end{proof}

The next result is a considerable
generalization of Theorem~\ref{thm02}
that encompasses Theorem~\ref{thmA} 
from the introduction.
Its proof requires
more notation.  

\begin{para} \label{paranew8}
Let $A$ be a normal domain and 
$m,n,r$ nonnegative integers such that
$r<\min\{m,n\}$.  Set $R=R_{r+1}(A;m,n)$ and consider the 
commutative diagram of natural ring homomorphisms
\begin{equation} \label{diag02} \tag{\ref{paranew8}.1} 
\begin{split}\xymatrix{
 & A[\mathbf{X}] \ar@{->>}[rd]^{\vf'} & \\
A \ar[ru]^{\dot{\vf}} \ar[rr]^{\vf} & & R
} \end{split} \end{equation}
wherein $\vf$ and $\dot{\vf}$ are faithfully flat,
and $\vf'$ is surjective and Cohen-Macaulay
of grade
$d=mn-r(m+n-r)$ by~\cite[(5.18)]{bruns:dr}; see~\ref{paranew6} for terminology.
If $C$ is a semidualizing $A$-module, then the following $R$-module
is semidualizing by~\cite[(6.1)]{frankild:appx}
\[ \cbc{C}{\vf}=\ext_{A[\mathbf{X}]}^d(R,C\otimes_A A[\mathbf{X}]).\]
 \end{para}

\begin{thm} \label{thm03}
Let $A$ be a 
normal domain and 
$m,n,r$ nonnegative integers such that
$r<\min\{m,n\}$.  The ring $R=R_{r+1}(A;m,n)$
is $\s_0$-finite
if and 
only if $A$ is so,
and the ordering on $\s_0(R)$ is transitive if and only if 
the ordering on $\s_0(A)$ is so.  
More specifically, one has the following cases.
\begin{enumerate}[\quad\rm(a)]
\item \label{item30}
If $r=0$ or $m=n$, then 
$\s_0(\vf)$ is a 
perfectly order-respecting bijection
$$\s_0(\vf)\colon\s_0(A)\xra{\approx}\s_0(R). $$ 
\item \label{item31}
If $r>0$ and $m\neq n$, then the 
assignment
\[ ([C]_A,0)\mapsto [\cbc{C}{\vf}]_R
\qquad\qquad
([C]_A,1)\mapsto [C\otimes_A R]_R  \]
describes a perfectly order-respecting bijection
$$h\colon\s_0(A)\times\{0,1\}\xra{\approx}\s_0(R).$$
\end{enumerate}
\end{thm}

The proof of this result is rather long, so it is presented at 
the end of the section
in~\ref{para01}.  For now we focus on some consequences of the theorem.

\begin{para} \label{paranew11}
Continue with the notation of~\ref{paranew10}.
Let $\n$ be a prime ideal of $A$ and 
consider the prime ideals $\fN=(\n,\mathbf{X})A[\mathbf{X}]$ and 
$\m=(\n,\mathbf{x})R$.
Localizing and completing the diagram~\eqref{diag02} 
yield similar commutative diagrams
\[ \xymatrix{
 & A[\mathbf{X}]_{\fN} \ar@{->>}[rd]^{\vf'_{\m}} & & & 
\comp{A[\mathbf{X}]_{\fN}}
\ar@{->>}[rd]^{\comp{\vf'_{\m}}} \\
A_{\n} \ar[ru]^{\dot{\vf}_{\fN}} \ar[rr]^{\vf_{\m}} & & R_{\m}
& \comp{A}=\comp{A_{\n}} \ar[ru]^{\comp{\dot{\vf}_{\fN}}} 
\ar[rr]^{\comp{\vf}=\comp{\vf_{\m}}} & & \comp{R_{\m}}=\comp{R}
} \]
For semidualizing $A_{\n}$- and $\comp{A}$-modules $C_0$
and $C_1$, respectively, we set 
\begin{align*}
\cbc{C_0}{\vf_{\m}}
  &=\ext_{A[\mathbf{X}]_{\fN}}^d(R_{\m},C_0\otimes_{A_{\n}} 
    A[\mathbf{X}]_{\fN}) 
  & \text{(semidualizing for $R_{\m}$)} \\
\cbc{C_1}{\comp{\vf}}
  &=\ext_{\comp{A[\mathbf{X}]_{\fN}}}^d(\comp{R},C_1\otimes_{\comp{A}} 
     \comp{A[\mathbf{X}]_{\fN}})  
  & \text{(semidualizing for $\comp{R}$)}
\end{align*}
These local constructions are discussed
extensively
in~\cite[Section 6]{frankild:appx}.  
\end{para}

What follows is the localized version of 
Theorem~\ref{thm03}.

\begin{cor} \label{cor05}
Let $A=\coprod_{i\geq 0}A_i$ be a graded normal 
domain 
with $(A_0,\n_0)$ local, and set $\n=\n_0+\coprod_{i\geq 1}A_i$.  
\begin{enumerate}[\quad\rm(a)]
\item \label{item40}
If $r=0$ or $m=n$, then 
$\s_0(\vf_{\m})$ is a 
perfectly order-respecting bijection
$$\s_0(\vf_{\m})\colon\s_0(A_{\n})\xra{\approx}\s_0(R_{\m})$$ 
\item \label{item41}
If $r>0$ and $m\neq n$, then the assignment 
\[ ([C]_{A_{\n}},0)\mapsto[\cbc{C}{\vf_{\n}}]_{R_{\m}} \qquad\qquad
([C]_{A_{\n}},1)\mapsto[C\otimes_{A_{\n}} R_{\m}]_{R_{\m}} \]
describes a perfectly order-respecting bijection
\[ 
\s_0(A_{\n})\times\{0,1\}\xra{\approx}\s_0(R_{\m}). \]
\end{enumerate}
\end{cor}

\begin{proof}
The following diagrams (one for 
each of our cases) commute.
\[ \xymatrix{
\s_0(A) \ar[r]^{\approx} \ar[d]^{\s_0(\vf)}_{\approx} 
& \s_0(A_{\n}) \ar[d]^{\s_0(\vf_{\m})}  &  
\s_0(A)\times\{0,1\} \ar[r]^{\approx} \ar[d]^{h}_{\approx} 
& \s_0(A_{\n})\times\{0,1\} \ar[d]^{h_{\m}} \\
\s_0(R) \ar[r]^{\approx} & \s_0(R_{\m})  & 
\s_0(R) \ar[r]^{\approx} & \s_0(R_{\m}) 
} \]
The four horizontal maps are 
perfectly order-respecting bijections by Proposition~\ref{prop05}, 
as are two of the vertical ones 
by Theorem~\ref{thm03}.  
Thus, the two remaining maps are so as well.
\end{proof}

\begin{cor} \label{cor06}
Let $A=\coprod_{i\geq 0}A_i$ be a graded super-normal 
domain 
with $(A_0,\n_0)$ local and complete.  Set
$\n=\n_0+\coprod_{i\geq 1}A_i$ and let $\m,\comp{A},\comp{R}$ be as 
in~\ref{paranew11}.
\begin{enumerate}[\quad\rm(a)]
\item \label{item42}
If $r=0$ or $m=n$, then 
$\s_0(\comp{\vf})$ is a 
perfectly order-respecting bijection
$$\s_0(\comp{\vf})\colon\s_0(\comp{A})\xra{\approx}\s_0(\comp{R}).$$ 
\item \label{item43}
If $r>0$ and $m\neq n$, then the assignment
\[ ([C]_{\comp{A}},0)\mapsto[\cbc{C}{\comp{\vf}}]_{\comp{R}} 
\qquad\qquad
([C]_{\comp{A}},1)\mapsto[C\otimes_{\comp{A}} \comp{R}]_{\comp{R}} \]
describes a perfectly order-respecting bijection
$$ 
\s_0(\comp{A})\times\{0,1\}\xra{\approx}\s_0(\comp{R}).$$ 
\end{enumerate}
\end{cor}

\begin{proof}
The proof is almost identical to the previous 
one, using Corollary~\ref{cor03} in place of Proposition~\ref{prop05}.
One needs only note that, since $A$ is 
super-normal, the same is true of 
$R$ by~\cite[(5.17),(6.12)]{bruns:dr}.
\end{proof}

The next step is to iterate the previous three results.

\begin{cor} \label{cor07}
Let $A$ be a 
normal domain and $t$ a positive 
integer.  For $l=1,\ldots,t$ fix integers $r_l,m_l,n_l$ such that 
$0\leq r_l<\min\{m_l,n_l\}$ and let $\mathbf{X}_{l\ast\ast}=\{X_{lij}\}$ be an 
$m_l\times n_l$ matrix of variables.  
Let $\mathbf{X}$ denote the entire list of variables $X_{111},\ldots, 
X_{tm_tn_t}$ and set
\[ R=A[\mathbf{X}]/\textstyle\sum_{l=1}^t 
\displaystyle I_{r_l+1}(\mathbf{X}_{l\ast\ast})\]
with $\mathbf{x}$ the image in $R$ of the sequence $\mathbf{X}$.
Let $s$ be the number of indices $l$ 
such that $r_l>0$ and $m_l\neq n_l$.  
\begin{enumerate}[\quad\rm(a)]
\item \label{item44}
There is a perfectly order-respecting bijection
$$ 
\s_0(A)\times\{0,1\}^s 
\xra{\approx}\s_0(R).$$
\item \label{item45}
With $A,\n$ as in Corollary~\ref{cor05} and $\m=(\n,\mathbf{x})R$,
there is a perfectly order-respecting bijection
$$ 
\s_0(A_{\n})\times\{0,1\}^s 
\xra{\approx}\s_0(R_{\m})$$  
\item \label{item46}
With $A,\n$ as in Corollary~\ref{cor06} and $\m=(\n,x)R$,
let $\comp{A}$ 
and $\comp{R}$ denote the $\n$-adic and $\m$-adic completions of $A$ 
and $R$, respectively.
There is a perfectly order-respecting bijection
$$ 
\s_0(\comp{A})\times\{0,1\}^s 
\xra{\approx}\s_0(\comp{R}).$$
\end{enumerate}
\end{cor}

\begin{proof}
Write $R_0=A$ and for $l=1,\ldots t$ set
$R_l=R_{r_l+1}(R_{l-1};m_l,n_l)$.  Then $R_t\cong R$ and
part~\eqref{item44} is proved by induction on $t$ using 
Theorem~\ref{thm03}.  
Parts~\eqref{item45} and~\eqref{item46} now follow from 
Proposition~\ref{prop05} and Corollary~\ref{cor03}, respectively.
\end{proof}

Before continuing, we present some notation.

\begin{para} \label{paranew9}
Let $A$ be a ring and fix an integer $m\geq 1$ and an $A$-regular sequence
$\mathbf{y}=y_{1},\ldots,y_{q}\in A$.
Set $n=m+q-1$ and let $\mathbf{X}$ be an $m\times n$ matrix of variables.
The discussion before and after~\cite[(2.14)]{bruns:dr}
exhibits a surjection
$$A[\mathbf{X}]/I_m(\mathbf{X})\twoheadrightarrow A/(\mathbf{y})^m$$
whose kernel is generated by an $A[\mathbf{X}]/I_m(\mathbf{X})$-regular sequence.

For $l=1,\ldots,t$ fix integers $m_l,q_l\geq 1$ and a sequence
$\mathbf{y}_{l\ast}=y_{l1},\ldots,y_{lq_l}\in A$.
Assume that the sequence $\mathbf{y}_{\ast\ast}$
is $A$-regular 
and set
$$ 
B(A,\mathbf{y},\mathbf{m},\mathbf{q})
=A/\textstyle\sum_{l=1}^{t}\displaystyle(\mathbf{y}_{l\ast})^{m_l}.$$
With $R$ as in Corollary~\ref{cor07},
tensoring
copies of the surjection from the previous paragraph provides a surjection
\begin{equation}\label{surj}\tag{\ref{paranew9}.1}
\tau\colon R \twoheadrightarrow B(A,\mathbf{y},\mathbf{m},\mathbf{q})
\end{equation}
whose kernel is generated by an $R$-regular sequence.

Let 
$B$ be a ring and $u$ a positive integer.
For $l=1,\ldots,u$ fix a positive integer $p_l$ and 
variables $\mathbf{Z}_{l\ast}=Z_{l1},\ldots,Z_{lp_l}$.  
We consider the ring
\[ 
S=S(B,\mathbf{p})=B[\mathbf{Z}_{1\ast}]/(\mathbf{Z}_{1\ast})^2\otimes_B
\cdots\otimes_BB[\mathbf{Z}_{u\ast}]/(\mathbf{Z}_{u\ast})^2 
\]
which can be thought of in several different ways.
Each ring $B[\mathbf{Z}_{l\ast}]/(\mathbf{Z}_{l\ast})^2$ 
is isomorphic to the 
trivial extension $B\ltimes B^{q_l}$, so there is an isomorphism
\[ S\cong(B\ltimes B^{q_1})\otimes_B\cdots\otimes_B(B\ltimes 
B^{q_u}) \]
Next, set $S_0=B$ and take successive trivial extensions  
$S_l=S_{l-1}\ltimes(S_{l-1})^{q_l}$.  
From the previous description, 
there is an 
isomorphism $S\cong S_u$.
Finally, let $\mathbf{Z}$ denote the full list of variables 
$\mathbf{Z}=Z_{11},\ldots,Z_{up_u}$
and let $\mathbf{z}$ denote the image in $S$ of the sequence 
$\mathbf{Z}$.  From the definition of $S$, one 
obtains the isomorphism
\[ S\cong B[\mathbf{Z}]/\textstyle\sum_{l=1}^u
\displaystyle(\mathbf{Z}_{l\ast})^2. 
\]
If $B$ is (complete) local  with maximal ideal $\fr$, then 
$S$ is (complete) local with maximal ideal $(\fr,\mathbf{z})S$.
\end{para}

The final result of this paper generalizes Corollary~\ref{cor08} and 
contains 
Theorems~\ref{thmB} and~\ref{thmC} from the introduction.

\begin{cor} \label{cor09}
With $A,\n$ as in Corollary~\ref{cor06},
let $t,u$ be  nonnegative integers.
For $l=1,\ldots,t$ fix positive integers $m_l,q_l$ and a sequence
$\mathbf{y}_{l\ast}=y_{l1},\ldots,y_{lq_l}\in\n\comp{A}$,  
and let $s$ be the number of indices $l$ 
such that $m_l,q_l>1$.  
For $l=1,\ldots,u$ fix a positive integer $p_l$, and  
let
$r$ denote the number of indices $l$ such that $p_l> 1$.
\begin{enumerate}[\quad\rm(a)]
\item \label{item77}
Set $\comp{B}=B(\comp{A},\mathbf{y},\mathbf{m},\mathbf{q})$
and $\comp{S}=S(\comp{B},\mathbf{p})$.
If $\mathbf{y}$ is $\comp{A}$-regular, then there is a
perfectly order-respecting bijection
$$\s_0(\comp{A})\times\{0,1\}^{r+s} 
\xra{\approx}\s_0(\comp{S}).$$
Furthermore, if $J\subset\comp{S}$ is an ideal primary to the maximal ideal of 
$\comp{S}$ and $C$ is a semidualizing $\comp{S}$-module, then
$e(J,C)=e(J,\comp{S})$.
\item \label{item77b}
Assume that $\mathbf{y}$ is an $A_{\n}$-sequence in $\n A_{\n}$,
and set $B'=B(A_{\n},\mathbf{y},\mathbf{m},\mathbf{q})$
and $S'=S(B',\mathbf{p})$.
There are 
perfectly order-respecting bijections
$$\s_0(A_{\n})\times\{0,1\}^{r+s} 
\xra{\approx}\s_0(S')\xra{\approx}\s_0(\comp{S}).$$ 
Furthermore, if $J\subset S'$ is an ideal primary to the maximal ideal of 
$S'$ and $C$ is a semidualizing $S'$-module, then
$e(J,C)=e(J,S')$.
\item \label{item78}
Assume that $\mathbf{y}$ is an $A$-regular sequence in $\n$,
and set $B=B(A,\mathbf{y},\mathbf{m},\mathbf{q})$
and $S=S(B,\mathbf{p})$.
There is a perfectly order-respecting  injection
\[ \s_0(A)\times\{0,1\}^{r+s} 
\hookrightarrow \s_0(S).
\]
\end{enumerate}
\end{cor}

\begin{proof}
We 
prove part~\eqref{item77};  argue similarly for
the other parts.
Let $\mathbf{Z}$ be as in~\ref{paranew9}.  
Then there are isomorphisms
$$\comp{S}\cong \comp{B}[\mathbf{Z}]/
\textstyle\sum_{l=1}^u\displaystyle(\mathbf{Z}_{l\ast})^2 
\cong 
\comp{A}[\mathbf{Z}]/(\textstyle\sum_{l=1}^{t}
\displaystyle(\mathbf{y}_{l\ast})^{k_l}
+\textstyle\sum_{l=1}^u\displaystyle(\mathbf{Z}_{l\ast})^2). 
$$
By Corollary~\ref{cor02}\eqref{item17}, the natural map
$\s_0(\comp{A})\to\s_0(\comp{A}[\mathbf{Z}])$ 
is a perfectly order-respecting bijection.
Pass to the ring $\comp{A}[\mathbf{Z}]$ and use the fact
that $\mathbf{Z}$ is $\comp{A}[\mathbf{Z}]/(\mathbf{y})$-regular,
to reduce 
to the case $u=0=r$, that is, $\comp{S}=\comp{B}$.

For $l=1,\ldots,t$ set $r_l=m_l-1$ and $n_l=m_l+q_l-1$, 
and let $R,\m,\comp{R}$ be as in 
Corollary~\ref{cor07}.
The surjection~\eqref{surj} 
completes to a surjection $\comp{\tau}\colon \comp{R}\to \comp{S}$
whose kernel is generated by 
a $\comp{R}$-regular sequence.  
The perfectly order-respecting bijections
in the next sequence are in
Corollary~\ref{cor07}\eqref{item46} and~\cite[(4.2)]{frankild:sdcms}
respectively
$$\s_0(\comp{A})\times\{0,1\}^{r+s} 
\xra{\approx}\s_0(\comp{R})\xra{\approx}\s_0(\comp{S}).$$
The equality of  multiplicities follows from Corollary~\ref{cor08}\eqref{item72a}.
\end{proof}

\begin{para} \label{paranew10}
To keep things tangible, we give an explicit description of the 
injection $$\s_0(A)\times\{0,1\}^{r+s} 
\hookrightarrow \s_0(S)$$ 
from the previous corollary.  (The two bijections are 
described analogously.)  
For $l=t+1,\ldots,t+u$ set $m_l=2$ and $q_l=p_{l-t}$.
Set $S_0=A$.  For $l=1,\ldots,t$
take quotients
$S_l=S_{l-1}/(\mathbf{y}_{l\ast})^{m_l}$, 
and for $l=t+1,\ldots,t+u$ 
take 
trivial extensions  
$S_l=S_{l-1}\ltimes(S_{l-1})^{p_l}$, so that $S\cong S_{t+u}$.  
Each 
homomorphism $\vf_{l-1}\colon S_{l-1}\to S_l$ 
induces an injective map:
\begin{enumerate}[\quad\rm(a)]
\item \label{item47}
If $q_l=1$ or $m_l=1$, 
then set $f_{l-1}=\s_0(\vf_{l-1})\colon \s_0(S_{l-1})\to 
\s_0(S_l)$;
\item \label{item48}
If $m_l,q_l>1$, then let $f_{l-1}\colon \s_0(S_{l-1})\times\{0,1\}\to 
\s_0(S_{l})$ be given by
\begin{align*}
([C]_{S_{l-1}},0)&\mapsto 
\begin{cases}
[\ext^{q_l}_{S_{l-1}}(S_l,C)]_{S_l} & \text{if $l\leq t$} \\
[\Hom_{S_{l-1}}(S_l,C)]_{S_l} & \text{if $l>t$}
\end{cases}
\\
([C]_{S_{l-1}},1)&\mapsto [C\otimes_{S_{l-1}}S_l]_{S_l}. 
\end{align*}
\end{enumerate}
The desired inclusion is exactly the composition
$f_{t+u-1}\cdots f_0$.
\end{para}

The calculations of this section motivate a
refinement of Question~\ref{q01}.

\begin{question} \label{q02}
If $R$ is a local ring, must the 
cardinalities of the sets $\s_0(R)$ and $\s(R)$ 
be powers of 2?
\end{question}

Paragraph~\ref{para02}
explains the need for the 
``local'' 
hypothesis.
Beyond the results of this 
section, evidence justifying this question
can be found in~\cite[(3.4)]{frankild:sdcms}:  
If $R$ is a non-Gorenstein ring 
admitting a dualizing complex and $\s(R)$ is a finite set, then 
$\s(R)$ has even cardinality.

We conclude this section with the proof of Theorem~\ref{thm03}.

\begin{para} \label{para01} \emph{(Proof of Theorem~\ref{thm03}.)}
Let $x_{ij}$ denote the residue of $X_{ij}$ in $R$ and set 
\[ \Delta=\det\begin{pmatrix} X_{11} & \ldots & X_{1r} \\ \vdots & 
 & \vdots \\ X_{r1} & \ldots & X_{rr} \end{pmatrix}\in A[\mathbf{X}]
\qquad 
\delta 
=\det\begin{pmatrix} x_{11} & \ldots & x_{1r} \\ \vdots & 
 & \vdots \\ x_{r1} & \ldots & x_{rr} \end{pmatrix}\in R \]
noting $\delta=\vf'(\Delta)$.
Also, set $e=\dim R -\dim A$ and
note that~\cite[(5.18)]{bruns:dr} implies $e=(m+n-r)r$.
By~\cite[(6.4)]{bruns:dr}
there is a prime element
$\zeta\in A[T_1,\ldots,T_e]$ and
an
isomorphism $\epsilon$ as in the next display; 
the isomorphism $\tau$ is clear.
\[ \xymatrix{
R\otimes_{A[\mathbf{X}]}A[\mathbf{X}]_{\Delta} \ar[r]^-{\cong}_-{\tau} &
R_{\delta} & 
A[T_1,\ldots,T_e]_{\zeta} 
\ar[l]_-{\cong}^-{\epsilon} 
} \]
Furthermore, the natural map 
$\alpha\colon A\to A[T_1,\ldots,T_e]_{\zeta}$ is 
faithfully flat so $\s_0(\alpha)$ is bijective by 
Corollary~\ref{cor02}\eqref{item17a}.

Set $U=A\smallsetminus(0)$ and $F=U^{-1}A$.
Using Lemma~\ref{prop04}\eqref{item15}, the natural maps
$\beta\colon R\to R_{\delta}$ and $\gamma\colon R\to U^{-1}R$
along with $\epsilon\alpha$
yield a
commutative diagram
\begin{equation} \tag{\ref{para01}.1}\label{diag01}
\begin{split}
\xymatrix{
\s_0(R) \ar[r]^-f \ar@{^{(}->}[d] & 
    \s_0(R_{\delta})\times\s_0(U^{-1}R) \ar@{^{(}->}[d] & \ar[l]_-g  
    \s_0(A)\times\s_0(U^{-1}R) \ar@{^{(}->}[d] \\
\Cl(R) \ar[r]^-{\cong}_-{f'} & 
    \Cl(R_{\delta})\times\Cl(U^{-1}R) & \ar[l]_-{\cong}^{-g'} 
    \Cl(A)\times\Cl(U^{-1}R) 
} \end{split} \end{equation}
where the horizontal maps are given by
\begin{gather} \tag{\ref{para01}.2}\label{diag01f}
\begin{split}
\xymatrix{
[C] \ar@{|->}[r]^-f & ([C\otimes_R R_{\delta}\,],[C\otimes_R  U^{-1}R]) & & &
} \\
 \xymatrix{
& & ([C'\otimes_AR_{\delta}],[C'']) & \ar@{|->}[l]_-g ([C'],[C''])
}  \end{split} \end{gather}
and the vertical arrows are induced by the respective inclusions.
The maps $f'$ and $g'$ are bijective 
by~\cite[(8.3)]{bruns:dr} and~\cite[(7.3),(8.1)]{fossum:dcgkd}, respectively.  
In particular, the maps $f,g$ are injective, and $g$ is bijective by 
Corollary~\ref{cor02}\eqref{item17a}.

\eqref{item30}  Assuming that $r=0$ or $m=n$,
Theorem~\ref{thm02} implies that $\s_0(U^{-1}R)$ 
is trivial 
since 
$U^{-1}R\cong R_{r+1}(F;m,n)$.  
Thus, the top row of~\eqref{diag01} reduces to
\begin{equation} \label{why1} \tag{\ref{para01}.3}
\xymatrix{
\s_0(R) \ar@{^{(}->}[r]^-{\s_0(\beta)} & \s_0(R_{\delta})
& \ar[l]_-{\s_0(\epsilon\alpha)}^-{\approx} \s_0(A).
} \end{equation}
The functoriality of $\s_0(-)$ 
and the following commutative diagram of ring homomorphisms
\[ \xymatrix{
A \ar[r]^-{\alpha} \ar[d]^{\vf} & A[T_1,\ldots,T_e]_{\zeta} \ar[d]^{\epsilon} \\
R \ar[r]^-{\beta} & R_{\delta}
} \]
yield the equality
$\s_0(\beta)\s_0(\vf)=\s_0(\epsilon\alpha)$.  Since~\eqref{why1} 
shows that $\s_0(\epsilon\alpha)$ is 
bijective, 
it follows that $\s_0(\beta)$ is surjective.  
As noted above,
$\s_0(\beta)$ is also 
injective, so it follows that 
$\s_0(\vf)$ is bijective. 
That it is a perfectly order-respecting bijection follows 
from~\cite[(4.8)]{frankild:appx} 
since $\vf$ is faithfully flat.

\eqref{item31}  Assume now that $r>0$ and $m\neq n$.  The isomorphism
$U^{-1}R\cong R_{r+1}(F;m,n)$ in conjunction with
Theorem~\ref{thm02} yields a bijection
$i\colon\{0,1\}\xra{\approx}\s_0(U^{-1}R)$ given by 
$i(0)=[\omega_{U^{-1}R}]$ 
and 
$i(1)=[U^{-1}R]$ where $\omega_{U^{-1}R}$ is a dualizing
module for 
$U^{-1}R$.  Let $i'\colon\s_0(A)\times\{0,1\}\to\s_0(A)\times\s_0(U^{-1}R)$ 
be the induced bijection.

Recall that $h\colon\s_0(A)\times\{0,1\}\to\s_0(R)$ is defined as
$h([C]_A,0)=[\cbc{C}{\vf}]_R$ and $h([C]_A,1)= [C\otimes_A R]_R$.
Below we construct a bijection 
\[ j\colon\s_0(R_{\delta})\times\s_0(U^{-1}R)\to 
\s_0(R_{\delta})\times\s_0(U^{-1}R) \]
such that the following diagram commutes.
\[ \xymatrix{
\s_0(A)\times\{ 0,1\} \ar[r]^-{h} \ar[d]_{i'}^{\approx} 
& \s_0(R) \ar@{^{(}->}[r]^-{f} & 
\s_0(R_{\delta})\times\s_0(U^{-1}R) \ar[d]_j^{\approx} \\
\s_0(A)\times\s_0(U^{-1}R) \ar[rr]^-{g}_-{\approx} & &
\s_0(R_{\delta})\times\s_0(U^{-1}R)
} \]
Once this 
is done, 
a simple diagram chase provides the bijectivity of $h$.  

Localize the surjection $\vf'\colon A[\mathbf{X}]\to R$ by inverting $\Delta$ 
to obtain a surjection $\rho\colon A[\mathbf{X}]_{\Delta}\to 
R_{\delta}$.  
We claim that $\rho$ is Gorenstein;
see~\ref{paranew6} for terminology.  
In~\ref{paranew8}
it is observed that $\vf'$ is Cohen-Macaulay of grade $d$.
Hence,
the same is true of 
$\rho$.  
The diagram~\eqref{diag02} 
fits in the next commutative diagram of ring homomorphisms.
\[ \xymatrix{
 & A[\mathbf{X}] \ar[r]^-{\psi} \ar@{->>}[d]^{\vf'} & A[\mathbf{X}]_{\Delta} 
\ar@{->>}[d]^{\rho} \\
A \ar[ru]^{\dot{\vf}} \ar[r]^-{\vf} \ar[rd]_{\alpha} & R \ar[r]^-{\beta}  
& R_{\delta} \\
& A[T_1,\ldots,T_e]_{\zeta} \ar[ru]^{\epsilon}_{\cong}
} \]
The map $\alpha$ is faithfully flat.  Furthermore, for each 
prime ideal $\p$ of $A$, the fibre 
$\kappa(\p)\otimes_A A[T_1,\ldots,T_e]_{\zeta}\cong
\kappa(\p)[T_1,\ldots,T_e]_{\zeta}$ 
is Gorenstein.  That $\rho$ is Gorenstein now follows 
from Avramov and Foxby~\cite[(6.2),(6.3)]{avramov:lgh}.

Set 
$\omega_{\rho}
=\ext^d_{A[\mathbf{X}]_{\Delta}}(R_{\delta},A[\mathbf{X}]_{\Delta})$,
which is $R_{\delta}$-semidualizing by~\cite[(5.6.a)]{frankild:appx}.
Moreover, it is locally free of rank 1 by~\cite[(5.6.b)]{frankild:appx}. 
Setting
\[ 
\omega_{\rho}^{-1}=\Hom_{R_{\delta}}(\omega_{\rho},R_{\delta})\] 
the discussion in~\ref{para02}  yields an isomorphism
\begin{equation}
\omega_{\rho}\otimes_{R_{\delta}}\omega_{\rho}^{-1}\cong 
R_{\delta}. \label{eq05} \tag{\ref{para01}.4}
\end{equation}
We now define the aforementioned map
$j$ and demonstrate that it has the 
desired properties.  For each semidualizing $R_{\delta}$-module 
$C$, set
\begin{align*}
j([C],[U^{-1}R])
&=([C],[U^{-1}R]) & &&
j([C],[\omega_{U^{-1}R}])
&=([\omega_{\rho}^{-1}\otimes_{R_{\delta}}C],
[\omega_{U^{-1}R}]).
\end{align*}
It follows from the isomorphism~\eqref{eq05} that the assignment
\begin{align*}
([C],[U^{-1}R])
&\mapsto([C],[U^{-1}R]) &&&
([C],[\omega_{U^{-1}R}])
&\mapsto([\omega_{\rho}\otimes_{R_{\delta}}C],
[\omega_{U^{-1}R}])
\end{align*}
describes an inverse of $j$, so that $j$ is bijective.  
It remains 
only to show that $gi'=jfh$, so fix a semidualizing $A$-module $C$.
First, there are 
isomorphisms
$$ 
C\otimes_A U^{-1}R\cong 
(C\otimes_A U^{-1}A)\otimes_{U^{-1}A} U^{-1}R\cong
U^{-1}A\otimes_{U^{-1}A} U^{-1}R 
\cong U^{-1}R
$$
the first and third of which are standard, and the second of which is
due to the fact that $U^{-1}A$ is a field.  This yields equality (1)
in the following sequence
\begin{align*}
jfh([C],1)
  & = jf([C\otimes_A R]) \\
  & = j([C\otimes_A R\otimes_R R_{\delta})],
      [C\otimes_A R\otimes_R U^{-1}R]) \\
  & \stackrel{(1)}{=}
     j([C\otimes_A R_{\delta}] ,
      [U^{-1}R]) \\
  & = ([C\otimes_A R_{\delta}] ,
      [U^{-1}R]) \\
  & = g([C],[U^{-1}R])\\
  & = gi'([C],1)
\end{align*}
where each of the unmarked equalities follows either from a definition 
(e.g., \eqref{diag01f}) or 
by a standard isomorphism.  

To compute $jfh([C],0)$, we first describe some isomorphisms:
\begin{align*}
\cbc{C}{\vf}\otimes_R R_{\delta} 
  & =
    \ext^d_{A[\mathbf{X}]}(R,C\otimes_A A[\mathbf{X}])\otimes_R R_{\delta} \\
  & \stackrel{(2)}{\cong}
    \ext^d_{A[\mathbf{X}]}(R,C\otimes_A A[\mathbf{X}])\otimes_R 
    (R\otimes_{A[\mathbf{X}]}A[\mathbf{X}]_{\Delta}) \\
  & \cong \ext^d_{A[\mathbf{X}]}(R,C\otimes_A A[\mathbf{X}]) 
    \otimes_{A[\mathbf{X}]}A[\mathbf{X}]_{\Delta}  \\
  & \stackrel{(3)}{\cong}
    \ext^d_{A[\mathbf{X}]_{\Delta}}
    (R\otimes_{A[\mathbf{X}]}A[\mathbf{X}]_{\Delta},
    C\otimes_A A[\mathbf{X}]\otimes_{A[\mathbf{X}]}A[\mathbf{X}]_{\Delta})  \\
  & \stackrel{(4)}{\cong}
    \ext^d_{A[\mathbf{X}]_{\Delta}}(R_{\delta}, 
    C\otimes_A A[\mathbf{X}]_{\Delta}) \\
  & \stackrel{(5)}{\cong}
    \omega_{\rho}\otimes_{R_{\delta}}
    (C\otimes_A A[\mathbf{X}]_{\Delta}\otimes_{A[\mathbf{X}]_{\Delta}}
    R_{\delta}) \\
  & \cong \omega_{\rho}\otimes_{R_{\delta}}
    (C\otimes_A R_{\delta}) 
\end{align*}
Each of the unmarked isomorphisms is either by definition or 
standard.  Isomorphisms (2) and (4) are via the isomorphism $\tau$, 
whereas (3) is from the flatness of $\psi$. For (5) use the equality
$\pd_{A[\mathbf{X}]_{\Delta}}(R_{\delta})=d$ to apply~\cite[(1.7.b)]{frankild:appx} 
and the definition of 
$\omega_{\rho}$.  Similar explanations yield all but one of the following isomorphisms. 
\begin{align*}
\cbc{C}{\vf}\otimes_R U^{-1}R
  & =\ext^d_{A[\mathbf{X}]}(R,C\otimes_A A[\mathbf{X}])\otimes_R U^{-1}R \\
  & \cong \ext^d_{A[\mathbf{X}]}(R,C\otimes_A A[\mathbf{X}])\otimes_R 
    (R\otimes_{A[\mathbf{X}]}U^{-1}A[\mathbf{X}]) \\
  & \cong \ext^d_{A[\mathbf{X}]}
    (R,C\otimes_A A[\mathbf{X}])\otimes_{A[\mathbf{X}]}U^{-1}A[\mathbf{X}] \\
  & \cong \ext^d_{U^{-1}A[\mathbf{X}]}
    (R\otimes_{A[\mathbf{X}]}U^{-1}\!A[\mathbf{X}],
    C\otimes_A A[\mathbf{X}]\otimes_{A[\mathbf{X}]}U^{-1}\!A[\mathbf{X}]) \\
  & \cong 
    \ext^d_{U^{-1}A[\mathbf{X}]}(U^{-1}R, C\otimes_A U^{-1}A[\mathbf{X}]) \\
  & \cong \ext^d_{U^{-1}A[\mathbf{X}]}(U^{-1}R,
    (C\otimes_A U^{-1}A)\otimes_{U^{-1}A} U^{-1}A[\mathbf{X}]) \\
  & \cong \ext^d_{U^{-1}A[\mathbf{X}]}(U^{-1}R,
    U^{-1}A\otimes_{U^{-1}A} U^{-1}A[\mathbf{X}]) \\
  & \cong \ext^d_{U^{-1}A[\mathbf{X}]}(U^{-1}R, U^{-1}A[\mathbf{X}]) \\
  & \stackrel{(6)}{\cong} \omega_{U^{-1}R}
\end{align*}
For isomorphism (6),
the ring $U^{-1}A[\mathbf{X}]$ is regular and surjects
onto $U^{-1}R$ so that 
$\ext^d_{U^{-1}A[\mathbf{X}]}(U^{-1}R, U^{-1}A[\mathbf{X}])$ is 
a dualizing module for $U^{-1}R$, and is therefore isomorphic to 
$\omega_{U^{-1}R}$ since the dualizing module of $U^{-1}R$ is unique 
up to isomorphism.

The preceding isomorphisms yield equality (7) in the next computation
\begin{align*}
jfh([C],0) 
  & = jf([\cbc{C}{\vf}]) \\
  & = j([\cbc{C}{\vf}\otimes_R R_{\delta}],
      [\cbc{C}{\vf}\otimes_R U^{-1}R]) \\
  & \stackrel{(7)}{=} 
    j([\omega_{\rho}\otimes_{R_{\delta}}
    (C\otimes_A R_{\delta})],
    [\omega_{U^{-1}R}]) \\
  & = ([\omega_{\rho}^{-1}\otimes_{R_{\delta}}
    \omega_{\rho}\otimes_{R_{\delta}}
    (C\otimes_A R_{\delta})],
    [\omega_{U^{-1}R}]) \\
  & \stackrel{(8)}{=} ([C\otimes_A R_{\delta}],
    [\omega_{U^{-1}R}]) \\
  & = g([C],[\omega_{U^{-1}R}])\\
  & = gi'([C],0).
\end{align*}
while (8) is by~\eqref{eq05}, and the others are by 
definition; see, e.g., \eqref{diag01f}.

To complete the proof, we verify the behavior of the orderings.
One implication follows from~\cite[(4.6),(5.7),(5.12)]{frankild:appx}: 
If $[C]_A\tri [C']_A$ and $i\leq i'$, then
$h([C]_A,i)\tri h([C']_A,i')$.  For the converse, 
assume that $h([C]_A,i)\tri h([C']_A,i')$.  By way
of contradiction, suppose that
$i>i'$, that is, $i=1$ and $i'=0$.  Then our assumption is
$[C\otimes_A R]_R\tri [C'(\vf)]_R$.  
The computations above provide isomorphisms
$$ C\otimes_A R\otimes_RU^{-1}R\cong U^{-1}R \qquad
\text{and}\qquad
C'(\vf)\otimes_RU^{-1}R\cong \omega_{U^{-1}R} $$
so that  the order-respecting map $\s_0(\gamma)$
yields $[U^{-1}R]_{U^{-1}R}\tri [\omega_{U^{-1}R}]_{U^{-1}R}$,
a contradiction since $U^{-1}R$ is not Gorenstein.
Thus, we have $i\leq i'$.  
The final conclusion
$[C]_A\tri [C']_A$ follows from~\cite[(4.8),(5.8),(5.13)]{frankild:appx}. \qed 
\end{para}

\section*{Acknowledgments}

I am grateful to Luchezar Avramov, Lars Christensen,
Neil Epstein, Anders Frankild,  Srikanth Iyengar, 
Graham Leuschke, Paul Roberts,
and Roger Wiegand 
for stimulating conversations  about this research.
Many thanks also to the referee for his/her very thorough comments.
I am especially grateful to Phillip Griffith for his tireless support and
motivation.  His work on the divisor class group
and my many conversations with him inspired the research contained
in this paper.

\providecommand{\bysame}{\leavevmode\hbox to3em{\hrulefill}\thinspace}
\providecommand{\MR}{\relax\ifhmode\unskip\space\fi MR }
\providecommand{\MRhref}[2]{%
  \href{http://www.ams.org/mathscinet-getitem?mr=#1}{#2}
}
\providecommand{\href}[2]{#2}

\end{document}